\documentclass[12pt, leqno]{amsart}
\usepackage{amsmath}
\usepackage{amssymb}
\usepackage{algorithm}
\usepackage[noend]{algpseudocode}

\usepackage{bm}
\usepackage{graphicx}
\usepackage{caption}
\usepackage{subcaption}
\usepackage[hidelinks]{hyperref}
\usepackage{mathtools}
\usepackage{comment}
\usepackage{enumerate} 
\usepackage{parskip}
\usepackage{hyperref}
 \usepackage{booktabs}

\usepackage{xcolor}
\allowdisplaybreaks

\newcommand{\R}{\mathbb{R}}

\newcommand{\E}{\mathbb{E}}

\newcommand{\Emb}{{\mathbb{E}}}
\newcommand{\Rmb}{{\mathbb{R}}}
\newcommand{\half}{{\frac{1}{2}}}
\newcommand{\dd}{\mathrm{d}}

\newtheorem{theorem}{Theorem}[section]

\newtheorem{proposition}[theorem]{Proposition}
\theoremstyle{definition}

\theoremstyle{remark}
\newtheorem{remark}[theorem]{Remark}
\numberwithin{equation}{section}

\renewcommand{\qedsymbol}{$\blacksquare$}

\numberwithin{equation}{section}

\begin{document}

\title[]{Neural Schr\"{o}dinger Bridge via Yasue Action with Applications to Generative Modeling}
\thanks{This research was partially supported by NSF award 2111688.}
\thanks{H. Liu acknowledges support from the Simons Foundation through its Travel Support for Mathematicians program.}

\author{Yan-Han Chen$^{\ast}$}

\author{Abhishek Halder$^{\ast}$}
\thanks{$^{\ast}$Department of Aerospace Engineering, Iowa State University, Ames, IA 50011, USA}

\author{Hailiang Liu$^{\dagger}$}
\thanks{$^{\dagger}$Department of Mathematics, Iowa State University, Ames, IA 50011, USA}



\begin{abstract}
We propose a novel neural solver for inferring the continuous state evolution of the Schr\"{o}dinger bridge problem. The proposed solver builds on a saddle point formulation involving the Yasue action functional. This allows us to realize the continuous state evolution by learning the optimal velocity for an ordinary differential equation, as opposed to learning the optimal drift for a stochastic differential equation. The proposed solver only uses the source and the target samples, and not the source and target distributions per se, which makes it particularly suitable for high dimensional applications arising in generative modeling. We demonstrate the effectiveness of the proposed solver using several numerical experiments in low and high dimensions, using both real and synthetic datasets.   
\end{abstract}

\maketitle

\section{Introduction}
\label{intro}
The Schr\"{o}dinger bridge (SB) \cite{leonard2014} -- now standard in generative modeling \cite{de2021diffusion,shi2023diffusion,JMLR:v26:23-1605} and stochastic control \cite{dai1991stochastic,chen2021stochastic,caluya2021wasserstein,caluya2021reflected} -- originated from the quest of reformulating quantum mechanics as non-equilibrium statistical mechanics \cite{schrodinger1931umkehrung,schrodinger1932theorie,wakolbinger1990schrodinger}. Two important follow ups in the spirit of this original quest are Nelson's stochastic mechanics \cite{nelson1967dynamical} and Yasue's least action principle \cite{yasue1981stochastic} -- the latter generalized the former by introducing minimization of an action functional over diffusion paths. In \cite{yasue1981stochastic}, Yasue noted: ``Clearly this least action principle is strongly analogous to Hamilton’s one in classical mechanics''. These ideas spurred several developments at the intersection of quantum mechanics, least action principle and stochastic control \cite{guerra1983quantization,carlen1984conservative,zambrini1986variational,pavon1995hamilton,vladimirov2012minimum,leger2019geometric,yordanov2024derivation,11208710}.

Intuitively, the SB can be understood as a minimum energy stochastic interpolant, i.e., an optimal measure-valued evolution between two given distributions at two given times. In this work, we propose a \emph{new neural solver for computing the SB dynamics via Yasue action}.

Motivated by generative modeling applications, our solver considers that only the samples from the endpoint distributions, and not the corresponding distributions themselves, are available. An interesting feature of the proposed approach is that the resulting measure-valued dynamics is realized by computing a \emph{deterministic velocity field} associated with an ordinary differential equation (ODE). This contrasts with standard computational methods for the SB involving drift and diffusion associated with a stochastic differential equation (SDE).

Our specific contributions are the following.
\begin{itemize}
    \item We propose a neural solver to learn the continuous state evolution in Schr\"odinger bridge problem that only uses the source and the target samples (scattered point clouds in $\mathbb{R}^d$). The proposed solver avoids spatial discretization. The solver builds on a saddle point formulation, and we prove the consistency for the same.
    
    \item We demonstrate the effectiveness of the proposed solver through extensive 
    numerical experiments on synthetic and real datasets, in both low and high-dimensional settings.
    
    \item Including the classical energy functional, the proposed solver can apply to general Lagrangians. 
\end{itemize}

\textbf{Related works.} Several computational methods exist for solving the SB problem and its generalizations, which do not use neural networks \cite{kullback1968probability,pavon2021data,caluya2021wasserstein,pooladian2025plug}. However, for high dimensional applications such as those arising in generative modeling, neural solvers have become the standard over the recent years. 

To compute the SB for generative modeling, De Bortoli et al. \cite{de2021diffusion} used score-based diffusion to approximate the well-known iterative proportional fitting procedure \cite{fortet1940resolution,kullback1968probability}. With similar motivation, Chen et al. \cite{chen2022likelihood} used forward-backward SDEs to train a pair of neural networks with  log-likelihood objective. For image-to-image translation, Su et al. \cite{su2023dual} proposed concatenating source
to latent, and latent to target Schr\"{o}dinger Bridges. See also \cite{I2SBpaper,kim2024unpaired}. Neural SB solvers have also been proposed for colloidal self-assembly \cite{nodozi2023neural}, for sampling \cite{liu2026adjoint}, and for single-cell perturbation prediction \cite{chi2026departures}. These references are only representative, not exhaustive, as this literature continues to grow very fast.

The existing SB solvers--both neural and non-neural--learn \emph{either} the optimal control drift in SDE, \emph{or} the solution for a system of partial differential equations (PDEs) that arise from the first order conditions for optimality, \emph{or} the likelihood, \emph{or} the static entropic potentials. Different from these, the proposed neural SB solver learns a deterministic velocity field via a saddle point reformulation of the SB problem involving the Yasue action.

\textbf{Notations.} We use $\mathbb{E}$ for the expectation operator, $\|\cdot\|$ for the Euclidean $2$-norm, $I$ for identity matrix of appropriate size, and $X\sim\rho$ to denote that the random variable $X$ has density $\rho$. The symbols $\nabla,\nabla\cdot,\Delta$ denote the standard Euclidean gradient, divergence and Laplacian operators, respectively, all with respect to the underlying spatial variable. The symbol $\mathcal{N}(\mu,\Sigma)$ denotes a multivariate Gaussian density with mean vector $\mu$ and (positive definite) covariance matrix $\Sigma$.

\textbf{Organization.} The rest of this work is structured as follows. In Section \ref{sec:SBandYasueAction}, we provide the necessary background on the SB and the Yasue action, and clearly state the learning problem of interest. In Section \ref{sec:MethodAndAlgorithm}, we develop the proposed solution method based on a saddle point reformulation and score dynamics. There, we also explain how the proposed method can be naturally extended for solving the SB problem with more general cost. Section \ref{sec:Numerical} provides four numerical case studies using both synthetic and real data benchmarks. Section \ref{sec:ConcludingRemarks} concludes the work. Proofs and additional numerical experiment are collected in the Appendices.

\section{Schr\"{o}dinger Bridge Problem and Yasue Action}\label{sec:SBandYasueAction}

\subsection{Background}\label{subsec:background}
The SB problem concerns with learning the minimum energy stochastic evolution in continuous time\footnote{The choice of unit time interval $[0,1]$ is without loss of generality. Any fixed compact time interval can be used.} $t\in[0,1]$ between a source density $\rho_{\text{src}}$ and a target density $\rho_{\text{tg}}$, where both the source and the target densities are supported over (subsets of) $\Rmb^d$. 
This can be formulated as the stochastic optimal control problem:
\begin{subequations}\label{Sch_Lagr}
\begin{gather}
    \inf_{u\in \mathcal{U}}\E \int_0^1 \frac{1}{2} 
    \left\lVert u\bigl(X(t),t\bigl) \right\rVert^2 \dd t,\\
    \dd X(t) = u\left(X(t),t\right) \dd t+\sigma\, \dd B(t),\label{SDEcontrolledBM}\\
         X(0)\sim \rho_{\text{src}},\quad 
         X(1)\sim \rho_{\text{tg}},
\end{gather}
\end{subequations}
where $\mathbb{E}$ denotes the expectation with respect to the law of $X(t)$, the symbol $B(t)$ denotes the standard Brownian motion in $\mathbb{R}^{d}$, and $X(t)$ is the to-be-determined It\^{o} diffusion such that its law interpolates between $\rho_{\text{src}}$ and $\rho_{\text{tg}}$ over the fixed time horizon $[0,1]$. 

In \eqref{Sch_Lagr}, the scalar diffusion coefficient $\sigma>0$ is fixed, and $\mathcal{U}$ denotes the set of $\mathbb{R}^{d}$-valued feedback control policies with finite energy. Specifically,
\begin{align}
\mathcal{U}:=\bigg\{u:\mathbb{R}^{d}\times [0,1]\to\mathbb{R}^{d} \mid \E \int_0^1 \frac{1}{2} 
    \left\lVert u\bigl(X(t),t\bigl) \right\rVert^2 \dd t < \infty\bigg\}.
\label{FeasibleControlSet}    
\end{align}
The minimizing control policy in \eqref{Sch_Lagr} determines the controlled Brownian motion $X(t)$.

We think of $X(t)$ in \eqref{Sch_Lagr} as the (stochastic) state of a particle in $\mathbb{R}^{d}$ at time $t\in[0,1]$, and aim to learn the densities along its trajectory for $0\le t\le 1$. The It\^{o} SDE \eqref{SDEcontrolledBM} describes the microscopic evolution of this particle's state. Problem \eqref{Sch_Lagr} has the following equivalent macroscopic formulation: 
\begin{subequations}\label{Sch_Euler}
    \begin{gather}
    \inf_{u,\rho}\int_0^1\int_{\R^d}
    \frac{1}{2} \lVert u(x,t)\rVert^2
    \rho(x,t)\,\dd x \,\dd t,\label{Schh_Euler_obj}\\
    \frac{\partial\rho}{\partial t} = -\nabla\cdot\left(u\rho\right)+\frac{\sigma^2}{2} \Delta \rho\label{eqn:FP},\\
         \rho(\cdot,0)= \rho_{\text{src}}(\cdot),\,
         \rho(\cdot,1)= \rho_{\text{tg}}(\cdot),
    \end{gather}
\end{subequations}
where $X(t)\sim\rho(\cdot,t)$.
Indeed, Problem \eqref{Sch_Euler} and Problem \eqref{Sch_Lagr} determine the same optimal density flow \cite[Section~4.5]{chen2021stochastic}.

We assume that both $\rho_{\text{src}}$ and $\rho_{\text{tg}}$ have finite second moments\footnote{This assumption is valid for our case where finitely many samples from $\rho_{\text{src}}$ and $\rho_{\text{tg}}$ are available. This is because being compactly supported is sufficient for the second moment to be finite.}. Under this assumption, the existence-uniqueness of solution for Problem \eqref{Sch_Euler} (equivalently for Problem \eqref{Sch_Lagr}) is guaranteed \cite[Thm. 3.2]{jamison1974reciprocal}, \cite[Section 10]{wakolbinger1990schrodinger}. Furthermore, with this assumption, for a given control in \eqref{FeasibleControlSet}, any feasible $\rho(\cdot,t)$ for Problem \eqref{Sch_Euler} is guaranteed to be almost everywhere positive for all $t\in[0,1]$. 

\subsection{Yasue action}\label{subsec:YasueAction}
We define the \emph{current velocity field} as\footnote{This is well-defined because of positivity of $\rho$ mentioned earlier.} 
\begin{align}
    v(x,t) := u(x,t) - \frac{\sigma^2}{2}\nabla\log \rho(x,t).
\label{defCurrentvelocity}
\end{align} 
By direct computation, it is easy to verify that Problem \eqref{Sch_Euler} admits the following reformulation \cite{chen2016relation}:
\begin{subequations}\label{sch_fish}
    \begin{gather}
    \inf_{v,\rho}\int_0^1\int_{\R^d}
   \left[ \frac{1}{2} \lVert v(x,t)\rVert^2 + \frac{\sigma^4}{8}\left\lVert\nabla\log\rho(x,t)\right\rVert^2\right]\label{sch_fish:obj}
    \rho(x,t)\,\dd x\, \dd t,\\
    \frac{\partial\rho}{\partial t} = 
    -\nabla\cdot\left(v\rho\right),\label{eqn:continuity}\\
         \rho(\cdot,0)= \rho_{\text{src}}(\cdot), \quad \rho(\cdot,1)= \rho_{\text{tg}}(\cdot).
    \end{gather}
\end{subequations}
In particular, the term $\int_{\mathbb{R}^{d}}\lVert\nabla\log\rho(x,t)\rVert^2\rho(x,t)\dd x$ is the well-known \emph{Fisher information} functional \cite[p. 279]{villani2021topics}. The Fisher information regularized objective \eqref{sch_fish:obj} is the \emph{Yasue action} \cite[p. 683]{chen2016relation}. The (non-autonomous) vector field $\nabla\log \rho(x,t)$ is often referred to as the \emph{score} \cite{hyvarinen2005estimation}.

We clarify here that the minimizing pair $(v,\rho)$ in \eqref{sch_fish:obj} exists whenever the minimizing pair $(u,\rho)$ in SB \eqref{Sch_Euler} exists, e.g., under our standing assumption that $\rho_{\text{src}},\rho_{\text{tg}}$ have finite second moments. Then, the minimum values of \eqref{Schh_Euler_obj} and \eqref{sch_fish:obj} coincide and is finite. As a consequence, the minimizing $v$ in \eqref{sch_fish:obj} must be of finite energy, i.e., must belong to the set $\mathcal{U}$ in \eqref{FeasibleControlSet}. Thus, there is no loss of generality in restricting the feasible $v$ for problem 
\eqref{sch_fish} to the set $\mathcal{U}$.

Problem \eqref{sch_fish} has the following microscopic or particle-level equivalence:
\begin{subequations} \label{sch_fish_SDE}
    \begin{gather}
    \inf_{v}\int_0^1\mathbb{E}_{X(t)\sim\rho}
   \left[ \frac{1}{2} \lVert v\left(X(t),t\right)\rVert^2 + \frac{\sigma^4}{8}\left\lVert\nabla\log\rho\left(X(t),t\right)\right\rVert^2\right]
    \dd t,\\
    \dd X(t) = v\left(X(t),t\right)\dd t,\label{Xt_ODE}\\
         X(0)\sim \rho_{\text{src}},\quad X(1) \sim \rho_{\text{tg}}\label{Xt_term}.
    \end{gather}
\end{subequations}
Observe that while \eqref{eqn:FP} is an advection-diffusion (i.e., Fokker-Planck) PDE, \eqref{eqn:continuity} is only advection (i.e., Liouville) PDE. In other words, the reformulation \eqref{sch_fish} allows us to transform the particle kinematics description from an SDE \eqref{SDEcontrolledBM} to an ODE \eqref{Xt_ODE}. This kinematic simplification comes at the expense of the additive Fisher information regularization in the objective \eqref{sch_fish:obj}.

The structure of Problem \eqref{sch_fish_SDE} motivates us to generalize the saddle point formulation proposed in our prior work \cite{chen2026neural},
in which an optimal transport\footnote{Optimal transport is known \cite{mikami2004monge,leonard2014,chen2016relation} to be the zero noise limit of SB, i.e., corresponds to the $\sigma=0$ case in \eqref{Sch_Euler} or \eqref{sch_fish}. In that case, the $X(t)$ kinematics \eqref{SDEcontrolledBM} already reduces to an ODE.} problem was reformulated as a min-max problem, which led to the development of a neural solver. However, the developments in that work were not for the objective (\ref{sch_fish:obj}), and it is far from obvious how to follow similar footsteps in this Yasue action setting.

Still, as shown in the next section, by characterizing the score dynamics in terms of the current velocity field, we are able to recast the constrained variational Problem \eqref{sch_fish_SDE} into a saddle point problem. This in turn allows us to learn two adversarially parameterized functions from only the samples from the source and the target densities: $\rho_{\text{src}}$ and $\rho_{\text{tg}}$. 

Denoting the source samples as $Z$ and the target samples as $Y$, we can now crystallize the problem statement as follows.

\colorbox{gray!20!white}{\parbox{\textwidth}{\textbf{Problem statement.} Given $Z\sim\rho_{\text{src}}$ and $Y\sim\rho_{\text{tg}}$, learn the optimal current velocity field $v^{*}$ solving Problem \eqref{sch_fish_SDE}, and thereby infer $X(t)\sim\rho(\cdot,t)$, $t\in[0,1]$.}}

The problem statement is pertient to generative modeling, where the goal is to draw new samples from a desired distribution that is available only through samples. 
A common strategy is to learn to map a simple reference (i.e., source) distribution that is easy to sample, to the desired (i.e., target) distribution.
In our framework, this reference is $\rho_{\text{src}}$ and the desired is $\rho_{\text{tg}}$.

In the next Section \ref{sec:MethodAndAlgorithm}, we develop a computational method and algorithm to solve the above problem. 


\section{Method and Algorithm}\label{sec:MethodAndAlgorithm}
In Section \ref{subsec:saddlepoint}, we detail the proposed saddle point reformulation. Section \ref{subsec:ScoreDynamicsAndAlgorithm} builds the proposed algorithm based on score dynamics. In Section \ref{subsec:generalcost}, we explain how the proposed methodology can be generalized for SB problems with more general costs. 

\subsection{Saddle point reformulation}\label{subsec:saddlepoint}
Let $\phi\in C_b(\Rmb^d;\mathbb{R})$, the vector space of all continuous and bounded functions from $\mathbb{R}^{d}$ to $\mathbb{R}$. We start by noting that two random variables $X(1)$ and $Y$ have the same distribution (thus the same density $\rho_{\text{tg}}$) if and only if \cite[Thm. 18.1]{jacod2004probability}, \cite[Thm. 1.2]{billingsley1999convergence}
\begin{align}
\Emb\left[\phi\left(X(1)\right)\right]=\Emb\left[\phi\left(Y\right)\right] \quad\forall\phi\in C_b(\Rmb^d;\mathbb{R}).
\label{SameDensityTestFunctionCondition}
\end{align}
Considering \eqref{SameDensityTestFunctionCondition} as an equality constraint, we notice that the Problem \eqref{sch_fish_SDE} is equivalent to the following:
\begin{subequations}\label{pb:saddle_supinside}
    \begin{gather}
    \inf_{v}\Emb\int_0^1
   \left[ \half \lVert v\left(X(t),t\right)\rVert^2 + \frac{\sigma^4}{8}\left\lVert\nabla\log\rho\left(X(t),t\right)\right\rVert^2\right]
    \dd t\nonumber\\
    \hspace*{1.3in} + \sup_{\phi}\left(\Emb\left[ \phi\left(X(1)\right)\right]-\Emb\left[\phi\left(Y\right)\right]\right),\label{eqn:adversarial_supinside}\\
    \dd X(t) = v\left(X(t),t\right)\dd t,\\
         X(0)\sim \rho_{\text{src}}.\label{init_constraint_supinside}
    \end{gather}
\end{subequations}
To see the equivalence between \eqref{pb:saddle_supinside} and \eqref{sch_fish_SDE}, notice that if
the equality in \eqref{SameDensityTestFunctionCondition} does not hold, then the supremum in \eqref{eqn:adversarial_supinside} equals $+\infty$. 

On the other hand, if the equality in \eqref{SameDensityTestFunctionCondition} holds (i.e., $Y\sim\rho_{\text{tg}}$), then the supremum equals zero. Because the first summand in \eqref{eqn:adversarial_supinside} is independent of $\phi$, we can then equivalently write \eqref{pb:saddle_supinside} as
\begin{subequations}\label{pb:saddle}
    \begin{gather}
    \inf_{v}\sup_{\phi}\quad\Emb\int_0^1
   \left[ \half \lVert v\left(X(t),t\right)\rVert^2 + \frac{\sigma^4}{8}\left\lVert\nabla\log\rho\left(X(t),t\right)\right\rVert^2\right]
    \dd t\nonumber\label{pb:saddle_obj}\\
    \hspace*{2in} + \left(\Emb\left[ \phi\left(X(1)\right)\right]-\Emb\left[\phi\left(Y\right)\right]\right),\label{eqn:adversarial}\\
    \dd X(t) = v\left(X(t),t\right)\dd t,\\
         X(0)\sim \rho_{\text{src}}.\label{init_constraint}
    \end{gather}
\end{subequations}

In the following Proposition \ref{pr:consistent} (proof in Appendix \ref{pf:consistent}), we relate the solution for Problem \eqref{pb:saddle} with that of Problem \eqref{sch_fish_SDE}.
\begin{proposition} \label{pr:consistent}
If Problem \eqref{pb:saddle} admits a saddle point $(v^*,\phi^*)$, then $v^*$ solves Problem \eqref{sch_fish_SDE}.
\end{proposition} 
The form of (\ref{eqn:adversarial}) is reminiscent of the framework of Generative Adversarial Network (GAN) \cite{goodfellow2014generative}.
Thus motivated, we parameterize the actor $v$ and critic $\phi$ with two neural networks $v_\theta$ and $\phi_\omega$, respectively. 
As is well-known \cite{de2021approx}, neural networks well-approximate continuous functions.

In our case, $v_\theta$ further leads to a parameterization of the trajectory $X(t) = X_{\theta}(t)$ by discretizing the ODE (\ref{Xt_ODE}) as 
\begin{equation}\label{Xt_discretization}
    X(t_{j+1}) = X(t_{j}) +  \frac{v_{\theta}\left(X(t_{j}),t_j\right)}{N}, 
\end{equation}
where the time interval $[0,1]$ is partitioned equally into $N$ subintervals with time stamps $t_j = \frac{j}{N}$, $j\in\{0,1,\hdots,N\}$.
We assign $X(0) = Z\sim\rho_{\text{src}}$ in (\ref{Xt_discretization}) to satisfy the initial constraint (\ref{init_constraint}).

\subsection{Score dynamics and proposed algorithm}\label{subsec:ScoreDynamicsAndAlgorithm}
To develop a computational method for solving \eqref{pb:saddle}, we make the following key observation: for a current velocity $v$ satisfying \eqref{eqn:continuity}, the evolution of the score $\nabla \log\rho$ along a particle trajectory $X(t)$ can be fully characterized by an ODE stated next in Proposition \ref{pr:score_ODE} (proof in Appendix \ref{pf:score_evo}).
\begin{proposition} \label{pr:score_ODE}
Given an initial stochastic state $X(0)\sim\rho(\cdot,0)$ and its dynamics $\dot{X}(t) = v\left(X(t),t\right)$,
 the evolution of the score $s(x,t) :=\nabla \log\rho (x,t)$ along a particle trajectory $X(t)$ satisfies almost surely, the ODE
\begin{equation}\label{score_ODE}
      \frac{\dd}{\dd t}s\left(X(t),t\right) = -   \left[\nabla v \left(X(t),t\right)\right]^{\top}s(X(t),t)  -(\nabla\nabla\cdot v )(X(t),t),
\end{equation}
where the the $(i,j)$-th entry in the Jacobian matrix $\nabla v$ is
$\frac{\partial v_i}{\partial x_j}$. 
\end{proposition}
\begin{remark}
We clarify that the dynamics \eqref{score_ODE} is an ODE, not a PDE, as it describes how $s$ evolves along a trajectory $X(t)$.
\end{remark}

Using $v_{\theta}$ in (\ref{score_ODE}), we approximate the evolution of the score $s(X(t),t) = s_\theta(X(t),t)$ along the particle trajectory by
\begin{align}\label{st_discretization}
s\left(X(t_{j+1}),t_{j+1}\right) = s\left(X(t_{j}),t_{j}\right) &-\frac{1}{N}\left(\left[\nabla v_{\theta}\left(X(t_{j}),t_j\right)\right]^{\top} s\left(X(t_j),t_{j}\right)
    \right.\nonumber\\
&\qquad\qquad\left.\vphantom{\left[\nabla v_{\theta}\left(X(t_{j}),t_j\right)\right]^{\top}}+\nabla\nabla\cdot v_{\theta}\left(X(t_j),t_{j}\right)\right),
\end{align}
where $s\left(X(0),0\right) = \nabla\log \rho_{\text{src}}\left(X(0)\right) = \nabla\log \rho_{\text{src}}(Z)$.

In the case $\rho_{\text{src}}$ is known (e.g., standard Gaussian), the initial score is analytically available ($s\left(X(0),0\right) = -Z$ for standard Gaussian). 
However, when only samples $Z\sim\rho_{\text{src}}$ are available, we may pretrain another neural network $s_\psi$ to learn its score from samples; see e.g, \cite{song2020sliced,vincent2011connection}.

 \begin{algorithm}[t!]
\caption{Training Procedure for Solving Problem~\eqref{pb:sch_learning}.} 
\label{alg:sch_learning}
\begin{algorithmic}[1]
\Require samples $Z\sim\rho_{\text{src}}$ and $Y\sim\rho_{\text{tg}}$; noise level $\sigma$;
 time steps $N$;
  initial score $\nabla\log\rho_{\text{src}}(\cdot)$ when $\rho_{\text{src}}$ known, otherwise $s_\psi(\cdot)$ (pre-trained by denoising score matching)
\Statex \textbf{Networks.} $v_\theta:\mathbb R^d\times[0,1]\to\mathbb R^d$;\quad
  $\phi_\omega:\mathbb R^d\to\mathbb R$
\Statex \textbf{Subroutine} \textproc{Rollout}$(\theta,Z)$ 
\State $X^{(0)}\gets Z$;\quad $s^{(0)}\gets \nabla\log\rho_{\text{src}}(Z)$ or $s_\psi(Z)$ 
\For{$j=0,\dots,N-1$}
  \State $v^{(j)}\gets v_\theta(X^{(j)},t_j)$
  \State $a^{(j)}\gets(\nabla_x v_\theta)^{\!\top}s^{(j)}$
  \State $b^{(j)}\gets\nabla_x(\nabla_x\!\cdot\!v_\theta)$
  \State $X^{(j+1)}\gets X^{(j)}+\tfrac1N v^{(j)}$;\quad
         $s^{(j+1)}\gets s^{(j)}-\tfrac1N (a^{(j)}+b^{(j)})$ \Comment{(\ref{Xt_discretization}) and (\ref{st_discretization})}
\EndFor
\State \Return $\{X^{(j)},v^{(j)},s^{(j)}\}_{j=0}^{N}$  
\Statex \textbf{Training}
\For{$k=1,\dots,N_v$}
  \State Draw $\{Z_i\}_{i=1}^{n}\sim\rho_{\text{src}}$, $\{Y_i\}_{i=1}^{n}\sim\rho_{\text{tg}}$
  \Statex \quad\emph{(a) critic ascent --- $\sup_\omega$ in (\ref{sch_learning_obj})}
   \For{$\ell=1,\dots,N_\phi$}
    \State $\{X^{(N)}_{i}\}_{i=1}^n\gets$\textproc{Rollout}$(\theta,Z_i)$ \textbf{with no grad}
    \State $\displaystyle \mathcal L_\omega\gets -\Bigl(\tfrac1n\textstyle\sum_i\phi_\omega\left(X^{(N)}_i\right)
           -\tfrac1n\textstyle\sum_i\phi_\omega\left(Y_i\right)\Bigr)$
    \State Update $\omega$ by descending $\displaystyle \mathcal L_\omega$.
  \EndFor
  \Statex \quad\emph{(b) descent in $\theta$ --- $\inf_{\theta}$ in (\ref{sch_learning_obj})}
  \State  $\{X^{(j)}_i,v^{(j)}_i,s^{(j)}_i\}_{i,j}\gets$\textproc{Rollout}$(\theta,Z_i)$

  \State $\displaystyle \mathcal L_\theta \gets \tfrac{1}{nN}\textstyle\sum_j\textstyle\sum_i\left[
        \tfrac12\left\lVert v^{(j)}_i\right\rVert^2+
        \frac{\sigma^4}{8}\left\lVert s^{(j)}_i\right\rVert^2\right]+
        \tfrac1n\textstyle\sum_i\phi_{\omega}\left(X^{(N)}_i\right)$
  \State Update $\theta$ by descending $\displaystyle \mathcal L_\theta$.
\EndFor
\State \Return $v_\theta$, $\phi_\omega$
\end{algorithmic}
\end{algorithm}

Finally, with neural parameterized $v_\theta$, $\phi_\omega$, and approximations (\ref{Xt_discretization}), (\ref{st_discretization}), we arrive at the following learning formulation for Problem \eqref{pb:saddle} (thereby for Problem \eqref{sch_fish_SDE} via Proposition \ref{pr:consistent}):
\begin{subequations}\label{pb:sch_learning}
    \begin{gather}
        \min_{\theta}\max_{\omega}\quad
    \Emb\int_0^1\left[\frac{1}{2}
     \left\lVert v_{\theta} \left(X(t),t\right)\right\rVert^2+ \frac{\sigma^4}{8}\left\lVert s_{\theta}\left(X(t),t\right)\right\rVert^2\right]\dd t\nonumber\\
     \qquad\qquad\qquad\qquad+\left(
\Emb\left[\phi_{\omega}\left(X(1)\right)\right]-
\Emb\left[\phi_{\omega}(Y)\right]\right),\label{sch_learning_obj}  \\
\dot{X}(t) = v_{\theta}\left(X(t),t\right)\label{sch_learning_ODE} ,\\
\dot{s}_{\theta}\left(X(t),t\right) = -\left[\nabla v_{\theta}\left(X(t),t\right)\right]^{\top} 
s_{\theta}\left(X(t),t\right)  -(\nabla\nabla\cdot v_{\theta} ) \left(X(t),t\right),
\end{gather}
\end{subequations}
where $ X(0) = Z$ and $s_{\theta}\left(X(0),0\right)=\nabla\log\rho_{\text{src}}(Z)$.
The proposed algorithm for solving this min-max problem is described in
Algorithm \ref{alg:sch_learning}.

For the application to generative modeling, we draw a new sample from $\rho_{\text{tg}}$ by mapping a new sample $Z\sim\rho_{\text{src}}$ to $X(1)$ with learned velocity $v_{\theta^*}$ according to \eqref{sch_learning_ODE}.
This generating process requires neither a stochastic differential equation nor a score network.

\begin{remark}
Implementing Algorithm \ref{alg:sch_learning} does not require the sample sets drawn from $\rho_{\text{src}}$ and $\rho_{\text{tg}}$ to have the same size, nor does it require the samples to be paired. The reason is that the two expectations in the second term in (\ref{sch_learning_obj}) can be computed separately.
\end{remark}
\subsection{Extension to general cost}\label{subsec:generalcost}
Our method extends naturally to SB problems with general cost $L(x,u,t)$, i.e., to problems of the form
\begin{equation}\label{Sch_general}
\begin{gathered}
    \inf_{u\in \mathcal{U}}\quad\E \int_0^1 L\left(X(t),u\left(X(t),t\right),t\right) \dd t,\\
    \dd X(t) = u\left(X(t),t\right) \dd t+\sigma\, \dd B_t,\\
         X(0)\sim \rho_{\text{src}},\, 
         X(1)\sim \rho_{\text{tg}}.
\end{gathered}
\end{equation}
Notice that the choice $L(x,u,t) = \frac{1}{2} \lVert u \rVert^2 $ reduces to the classical SB problem \eqref{Sch_Lagr}. Problems of the form \eqref{Sch_general} appear naturally in stochastic control \cite{11239424,11735073,teter2024weyl,teter2025probabilistic} and learning \cite{koshizuka2023neural,liu2024generalized,behjoo2025harmonic,buzun2026hota,NEURIPS2024_78fa1d56}.

Following what we did in transcribing Problem \eqref{Sch_Euler} to
Problem \eqref{sch_fish_SDE} with the current velocity field $v(x,t)$, we may consider the following ``Yasue action reformulation" for Problem \eqref{Sch_general}:
\begin{equation}\label{Sch_general_fish}
    \begin{gathered}
         \inf_{v}\quad\Emb \int_0^1 L\Bigl(X(t),v\left(X(t),t\right)+\frac{\sigma^2}{2}\nabla\log\rho\left(X(t),t\right),t\Bigr)\dd t,\\
    \dd X(t) = v\left(X(t),t\right) \dd t,\\
         X(0)\sim \rho_{\text{src}},\,\, 
         X(1)\sim \rho_{\text{tg}}.
    \end{gathered}
\end{equation}
As before, Problems \eqref{Sch_general} and \eqref{Sch_general_fish} determine common optimal density flow. Furthermore, retracing the same steps as in Sections \ref{subsec:saddlepoint} and \ref{subsec:ScoreDynamicsAndAlgorithm}, the corresponding neural learning formulation for Problem \eqref{Sch_general_fish} is 
\begin{subequations}
\label{sch_learning_L}
    \begin{gather}
        \min_{\theta}\max_{\omega}
    \quad\Emb\int_0^1 L\Bigl(X(t),v_{\theta}\left(X(t),t\right)+\frac{\sigma^2}{2}s_{\theta}\left(X(t),t\right),t\Bigr)\dd t\nonumber\\
    \qquad\qquad\qquad\qquad+\left(
\Emb\left[\phi_{\omega}\left(X(1)\right)\right]-
\Emb\left[\phi_{\omega}(Y)\right]\right),  \\
\dot{X}(t) = v_{\theta}\left(X(t),t\right),\\
\dot{s}_{\theta}\left(X(t),t\right) = -\left[\nabla v_{\theta}\left(X(t),t\right)\right]^{\top} 
s_{\theta}\left(X(t),t\right)  -(\nabla\nabla\cdot v_{\theta} ) \left(X(t),t\right),
\end{gather}
\end{subequations}
where $ X(0) = Z$ and $s_{\theta}\left(X(0),0\right)=\nabla\log\rho_{\text{src}}(Z)$. It follows that the proposed Algorithm \ref{alg:sch_learning} can be applied mutatis mutandis to \eqref{sch_learning_L}.


\section{Numerical Experiments}\label{sec:Numerical}
In this Section, we evaluate our proposed Algorithm \ref{alg:sch_learning} on both synthetic and real, low and high dimensional datasets. 
Specifically, Section \ref{sec:MixG2D} illustrates the proposed method on a synthetic Gaussian source to Gaussian-mixture target case study. Section \ref{test:MNIST}, \ref{test:MNIST_Gau} and \ref{fashionMNIST} demonstrates the proposed method with MNIST~\cite{lecun1998gradient}, MNIST with Gaussian source, and Fashion-MNIST datasets~\cite{xiao2017fashion}, respectively. An additional numerical experiment with CelebA dataset~\cite{liu2015faceattributes} is reported in Appendix~\ref{app:CelebA}. In Section \ref{subsec:OTcomparison}, we compare the performance of the proposed neural SB solver with its zero noise ($\sigma=0$), i.e., optimal transport limit. Section \ref{subsec:CommentsNumerics} comments on the reported results.

Across all numerical experiments, we set the noise $\sigma = 1$.
Updating the score in~\eqref{st_discretization} involves two terms that are computationally expensive in high dimensions: $\nabla\cdot v_\theta$ and $\nabla_x v_\theta$.
For $\nabla\cdot v_\theta$, we follow \cite{grathwohl2018ffjord} and many other works to adopt the Skilling--Hutchinson trace estimator \cite{hutchinson1989stochastic,skilling1989eigenvalues}, which requires much lower computation for large $d$ at the price of small random noise.
For $\nabla_x v_\theta$, which is required only through its product with a vector, we circumvent computing the full Jacobian matrix by using a vector--Jacobian product (VJP).

The algorithm is implemented in PyTorch \cite{paszke2019pytorch}. 
Both neural networks
$v_\theta$ and $\phi_\omega$ are trained with the RMSprop optimizer \cite{tieleman2012lecture}, with gradients computed by reverse-mode automatic
differentiation in PyTorch's \texttt{autograd} engine \cite{paszke2017automatic}.
All experiments were run on the Nova high-performance computing cluster at Iowa State University using NVIDIA GPUs.
The details for the network architectures and training hyperparameters are given in Appendix~\ref{appdx:details_exp_hyper}.

\subsection{Gaussian mixture}\label{sec:MixG2D}
\begin{figure}
    \centering
\vspace{-4pt}
\includegraphics[width=\linewidth]{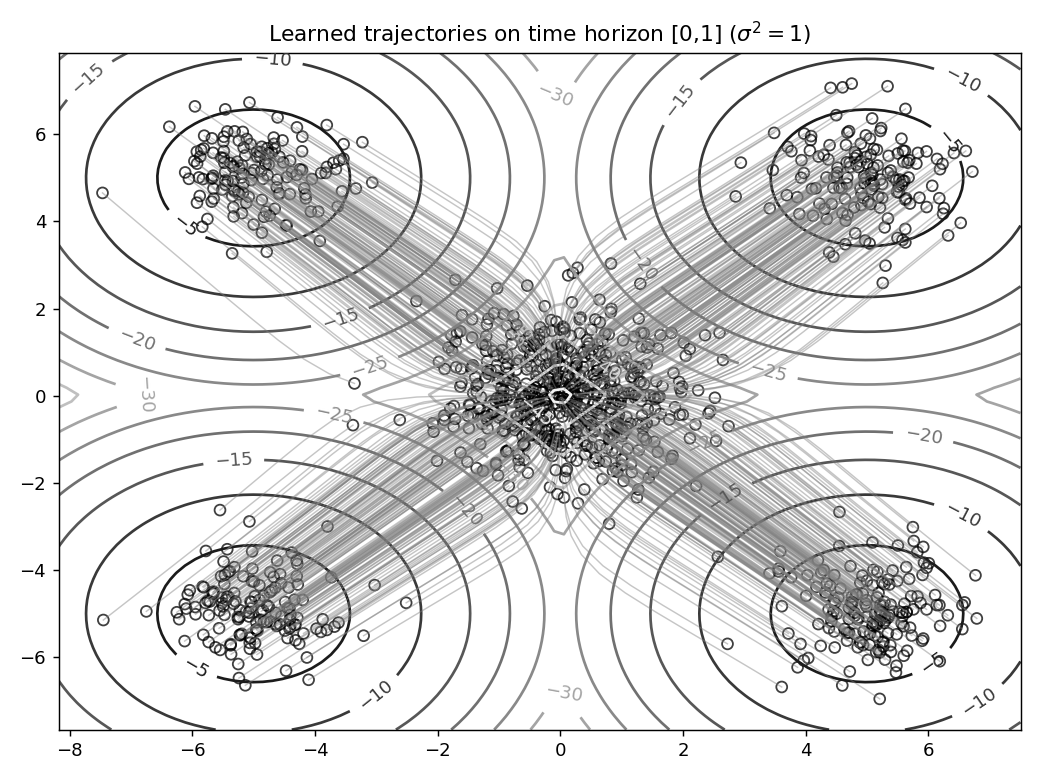}
    \captionsetup{aboveskip=-1pt} 
    \caption{Samples of the standard normal source $\rho_{\text{src}}$ (open circular markers at the center) are transported to the location of the four modes of the mixture of Gaussian $\rho_{\text{tg}}$. 
        Each gray line illustrates the learned trajectory obtained from our proposed Algorithm \ref{alg:sch_learning}. 
        The contour shows the 
        true log density of $\rho_{\text{tg}}$.
       }
    \label{fig:MixG_Xt}
\end{figure}

\begin{figure}
    \centering
\vspace{-1pt}
\includegraphics[width=\textwidth]{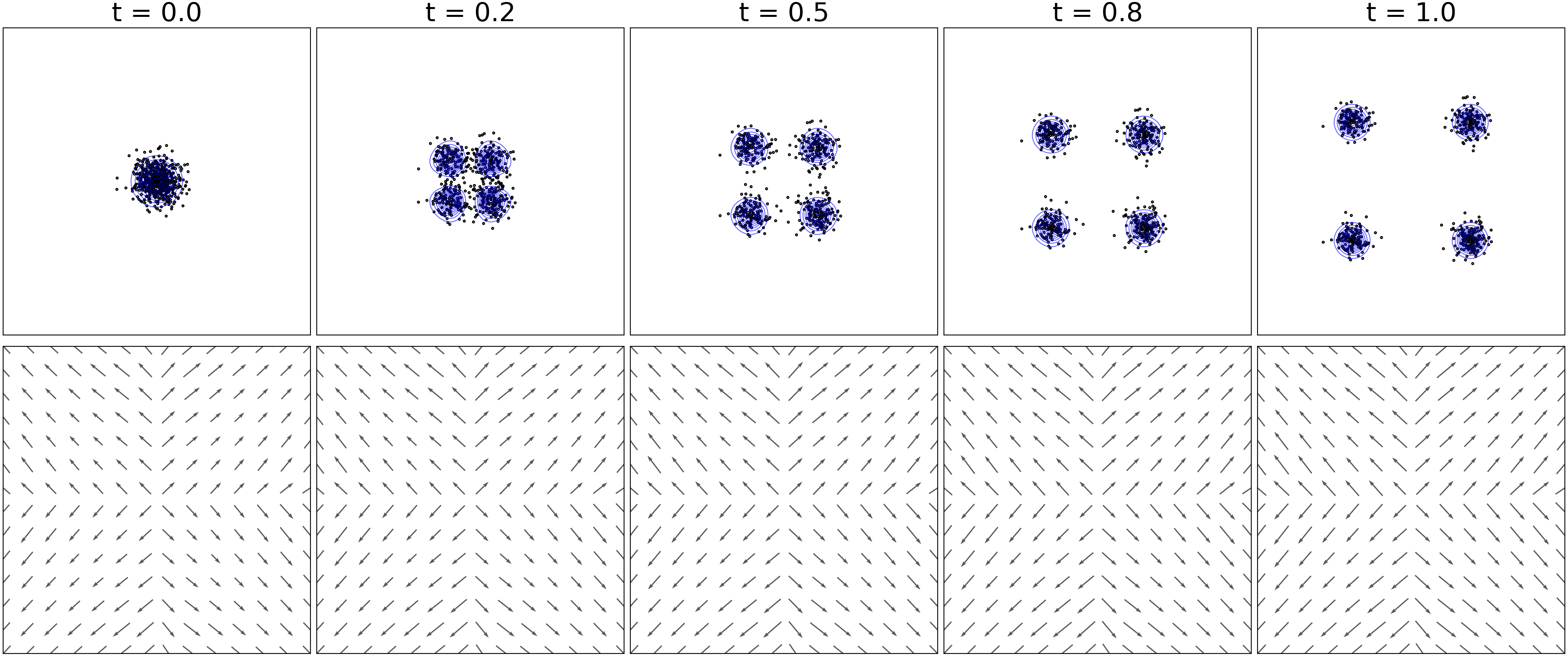}
    \caption{First row: the generated sample cloud from our method (open circular markers) aligns with the density contours obtained using the non-neural method in \cite{caluya2021wasserstein} and \cite{11239424} at selected time slices. 
    Second row: learned time-varying velocity field at selected time slices.}
    \label{fig:MixG_SBP}
\end{figure}

\begin{figure}
    \centering
    \vspace{-5pt}
    \includegraphics[width=\linewidth]{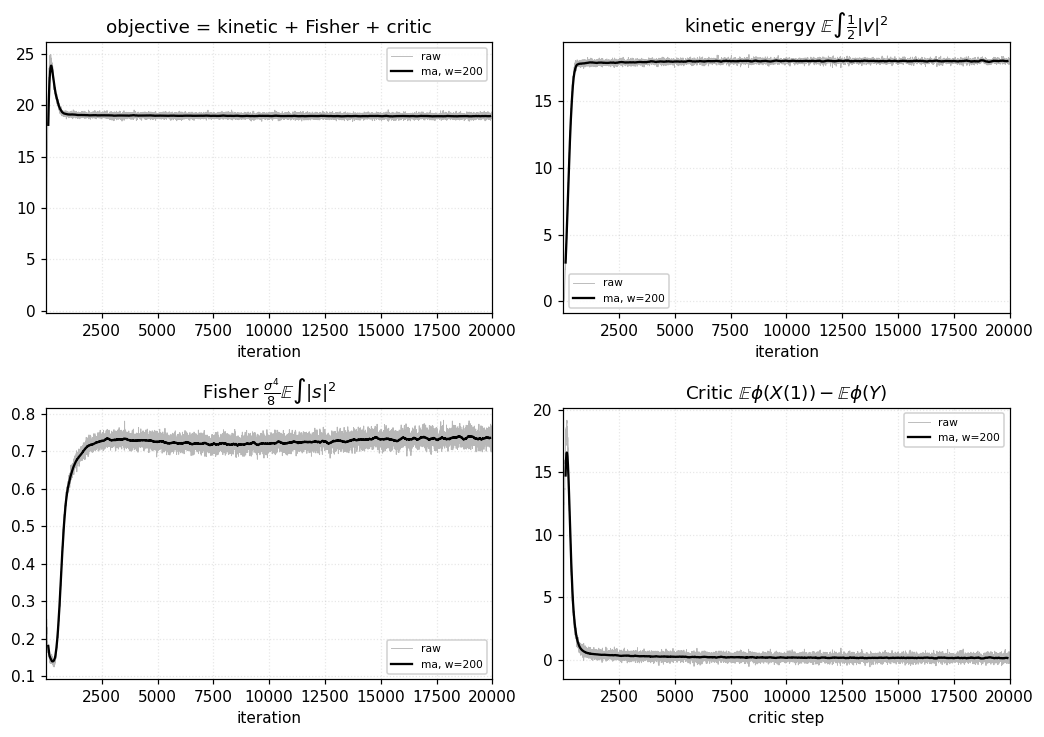}
    \captionsetup{aboveskip=0pt}
    \caption{
    The history of training process of Gaussian mixture example reporting the objective \eqref{sch_learning_obj} and the three terms in it. 
    The solid line is obtained by taking the moving average of the raw values in light gray, with window size 200 steps.}
   \label{fig:MixG_Xt_loss}
\end{figure}

In this synthetic experiment in $d=2$ dimensions, we consider the source density $\rho_{\text{src}}$ to be standard 2-dimensional Gaussian, and the target density $\rho_{\text{tg}}$ to be an equally weighted mixture of four Gaussian components. 
Specifically,
$$\rho_{\text{tg}}=\frac{1}{4}\sum_{i=1}^{4}\:\mathcal{N}(\mu_i,0.5 I),\quad \mu_i \in \{-5,5\}^2\quad\forall i=1,2,3,4.$$
Figure \ref{fig:MixG_Xt} shows the learned trajectories. 
To validate our neural algorithm, we compare the sample cloud generated by our method with the density contours computed by the non-neural method in \cite{caluya2021wasserstein} and \cite{11239424} at selected time slices. 
This comparison is reported in the first row of Figure \ref{fig:MixG_SBP}; the second row of that same figure shows the learned time-varying velocity field $v_{\theta^{*}}$ snapshots as quiver plots. Figure \ref{fig:MixG_Xt_loss} depicts the history of the neural network training for this experiment. We report the training history only for this example since this is representative in the sense the corresponding histories for the following experiments have similar trends.

\subsection{MNIST}\label{test:MNIST}
We next apply our algorithm to the MNIST handwritten-digit dataset \cite{lecun1998gradient}, which also provides a high-dimensional benchmark. Each digit is represented as a ($28 \times 28$) grayscale image and thus as a vector in $\mathbb R^{784}$. Specifically,  we  learn the current velocity field mapping from digit 9 (representing $\rho_{\text{src}}$) to digit 6 (representing $\rho_{\text{tg}}$). Since the true source score required in Algorithm \ref{alg:sch_learning} is unavailable in this setting, we pretrain the score function $s\left(X(0),0\right)$ using denoising score matching \cite{vincent2011connection}.

\begin{figure}[htbp]
     \centering
     \captionsetup[subfigure]{labelformat=empty} 
     
     \subfloat{%
         \includegraphics[width=0.495\textwidth]{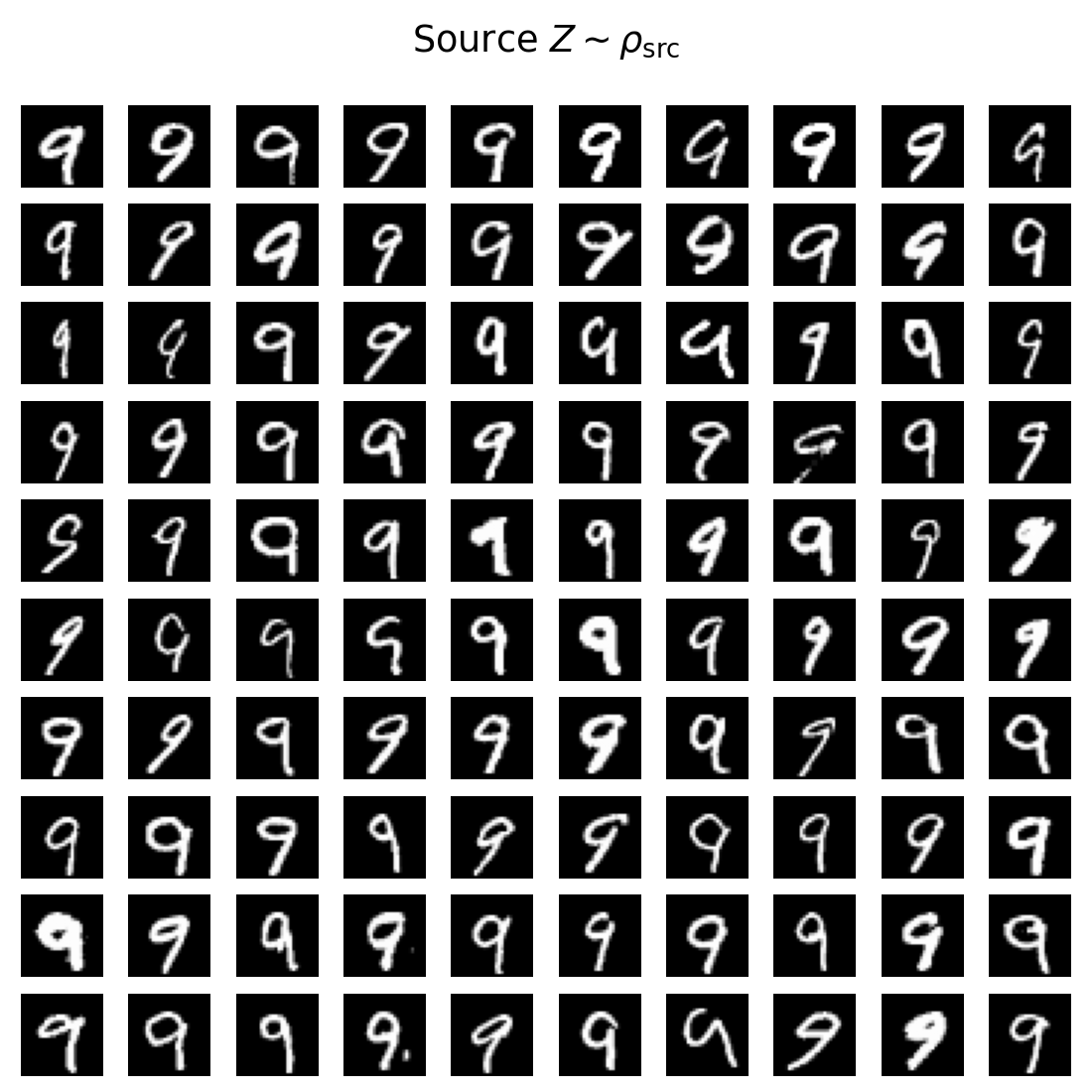}%
     }%
     \subfloat{%
         \includegraphics[width=0.495\textwidth]{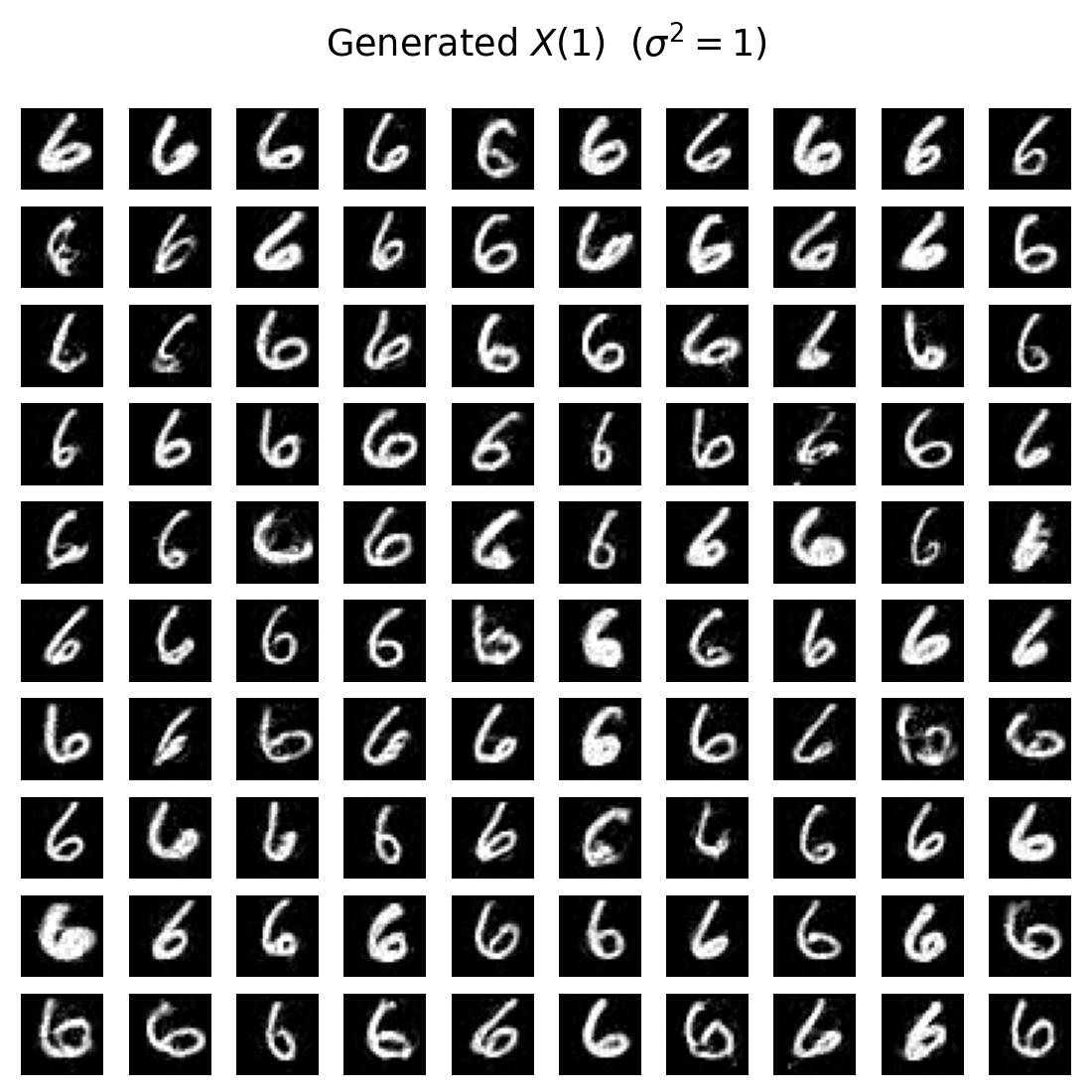}%
     }%
     \caption{
     Left panel: 100 samples from the source distribution $\rho_{\text{src}}$ (digit 9). 
     Right panel: transported samples from the left panel meant to match the target distribution $\rho_{\text{tg}}$ (digit 6).
     }
     \label{fig:MNIST_X1}
\end{figure}

\begin{figure}
    \centering
    \vspace{-5pt}
    \includegraphics[width=\linewidth]{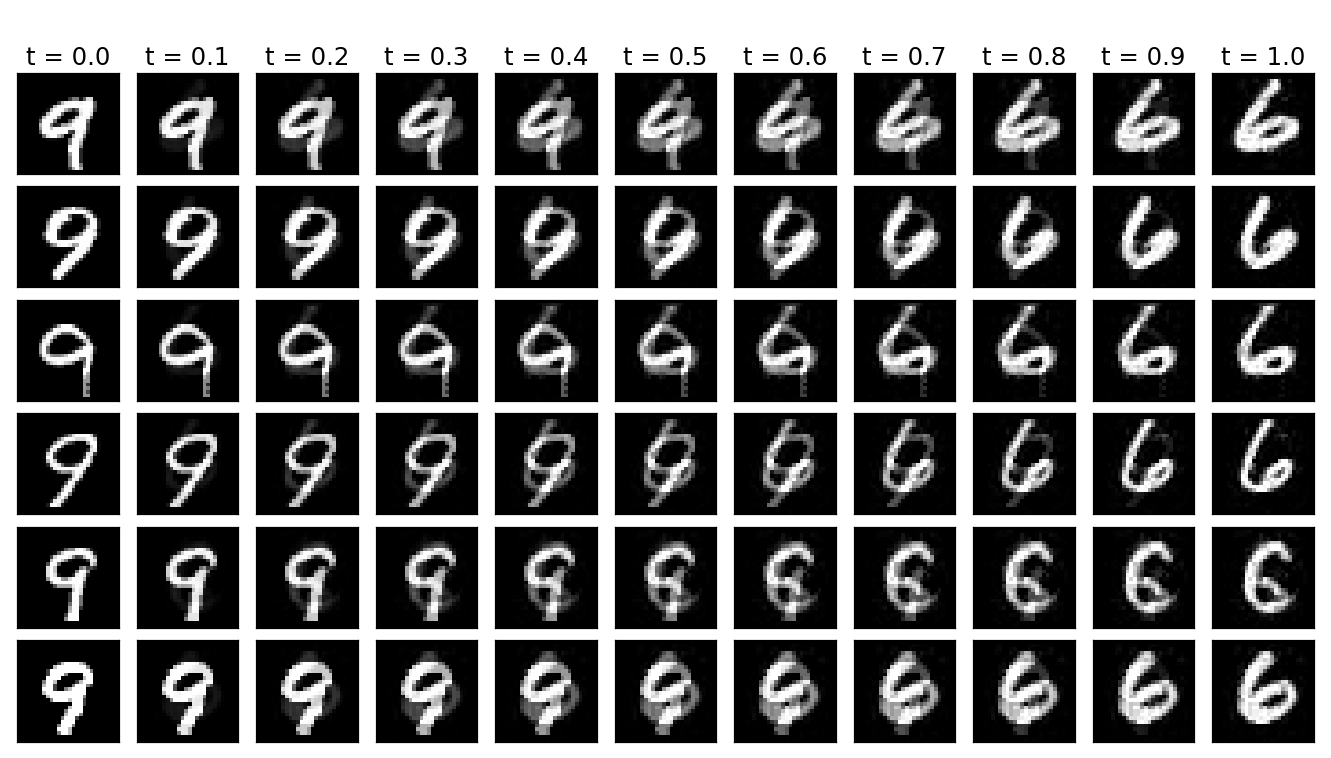}
    \captionsetup{aboveskip=0pt}
    \caption{
    Selected time slices of the transformation process $X(t)$.}
   \label{fig:MNIST_Xt}
\end{figure}

Figure \ref{fig:MNIST_X1} shows the generated digit (6) samples obtained from digit (9) samples using Algorithm \ref{alg:sch_learning}. Figure \ref{fig:MNIST_Xt} demonstrates the transformation induced by the learned velocity field $v_{\theta^{*}}$.


\subsection{MNIST with Gaussian source}\label{test:MNIST_Gau}
In this experiment, we again consider the MNIST dataset as in Section \ref{test:MNIST}, but now with a standard Gaussian $\rho_{\text{src}}$. The motivation behind this experiment is that transporting a Gaussian to digit $9$ is more challenging than transporting digit $6$ to digit $9$ (due to shape similarity). In this case, the source score is available analytically as 
$$
s(Z) = -Z, 
$$
and thus requires no approximation or pretraining.

\begin{figure}[htbp]
     \centering
     \captionsetup[subfigure]{labelformat=empty} 
     
     \subfloat{%
         \includegraphics[width=0.495\textwidth]{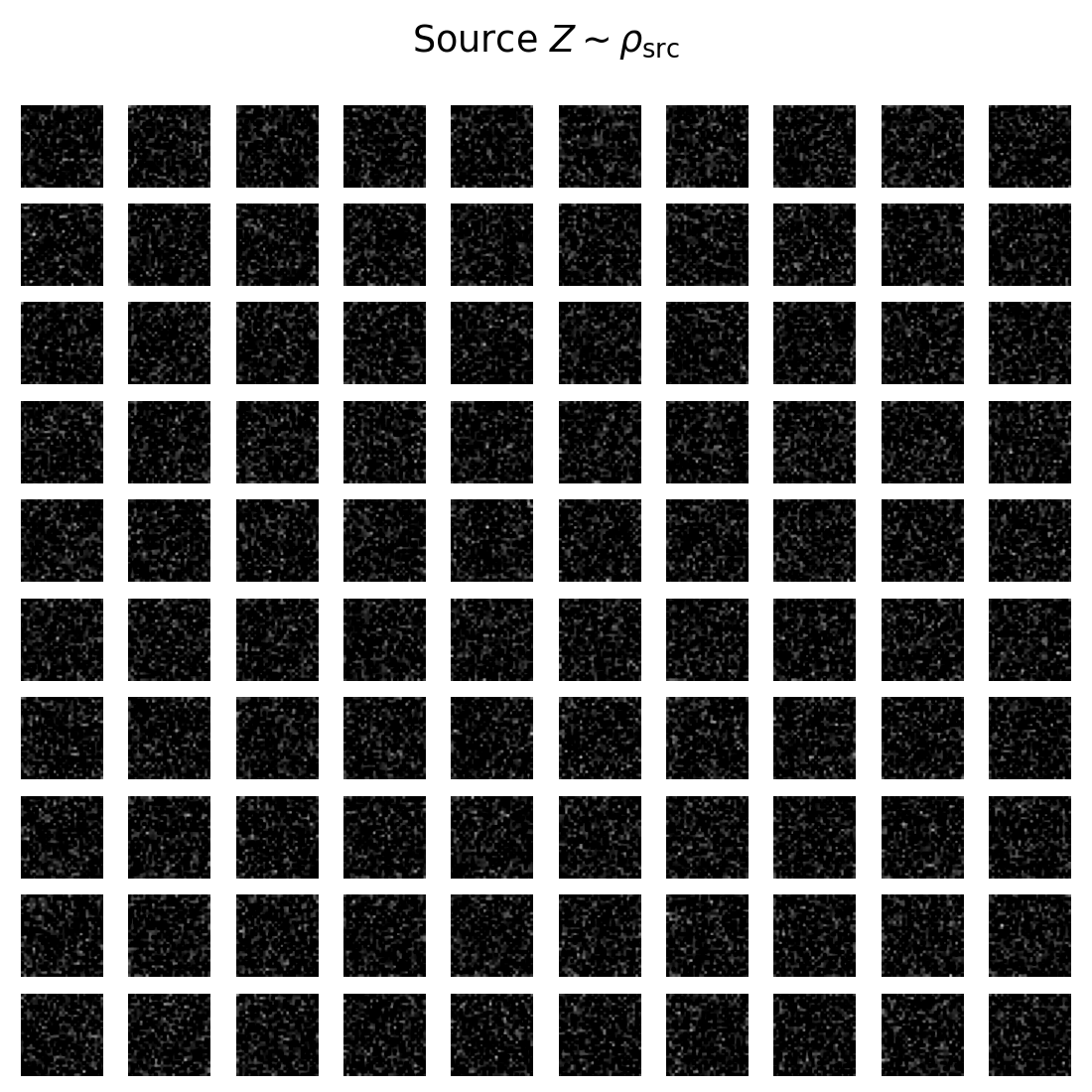}%
     }%
     \subfloat{%
         \includegraphics[width=0.495\textwidth]{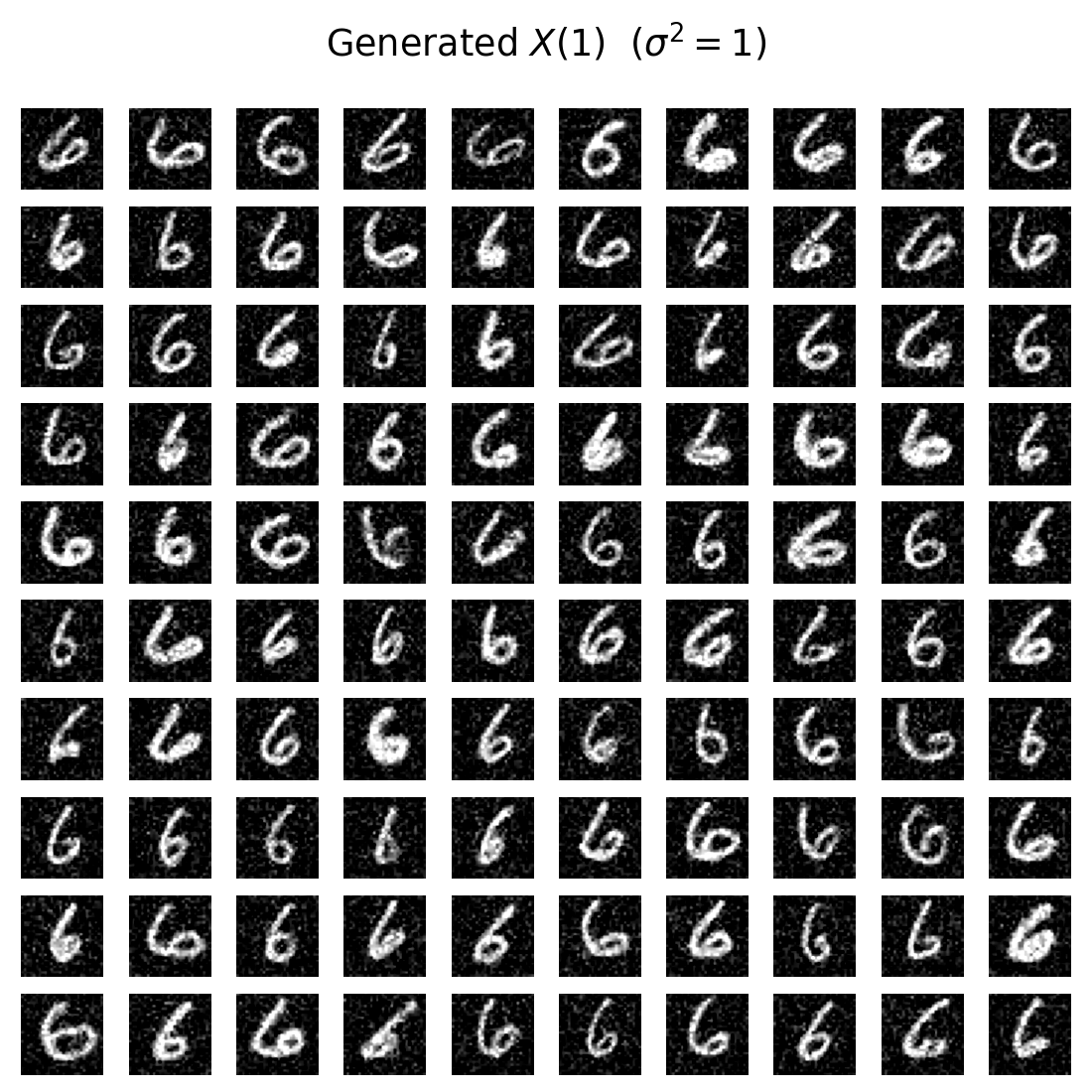}%
     }%
     \caption{
     Left panel: 100 samples from the source distribution $\rho_{\text{src}} = \mathcal{N}(0,I)$. 
     Right panel: transported samples from the left panel meant to match the target distribution $\rho_{\text{tg}}$ (digit 6).
     }
     \label{fig:MNISTG_X1}
\end{figure}

\begin{figure}
    \centering
    \vspace{-5pt}
    \includegraphics[width=\linewidth]{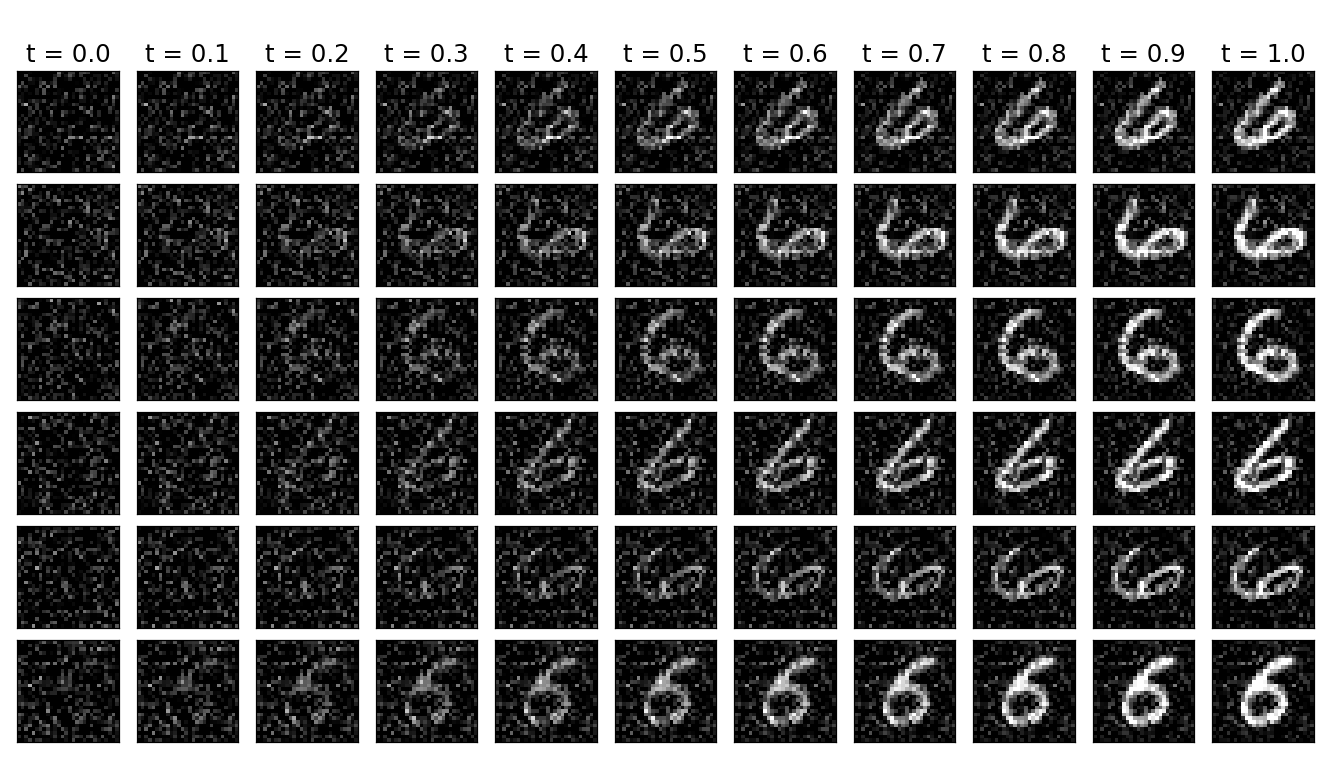}
    \captionsetup{aboveskip=0pt}
    \caption{
    Selected time slices of the transformation process $X(t)$.}
   \label{fig:MNISTG_Xt}
\end{figure}

Figure \ref{fig:MNISTG_X1} reports the generated digit (6) obtained by transporting the samples drawn from the standard Gaussian distribution via the learned velocity field $v_{\theta^{*}}$ from the proposed Algorithm \ref{alg:sch_learning}. Figure \ref{fig:MNISTG_Xt} demonstrates the transformation induced by $v_{\theta^{*}}$.

\subsection{Fashion MNIST}\label{fashionMNIST}
As a further real data experiment, we test our method on the Fashion-MNIST dataset \cite{xiao2017fashion}, which  consists of 70,000 grayscale images of size ($28\times28$).
We learn a transport map from the distribution of ``sneaker'' (representing $\rho_{\text{src}}$) to that of ``ankle boot'' (representing $\rho_{\text{tg}}$).  Since the true source score is unavailable, we again use denoising score matching to estimate it. 

\begin{figure}[htbp]
     \centering
     \captionsetup[subfigure]{labelformat=empty} 
     
     \subfloat{%
         \includegraphics[width=0.495\textwidth]{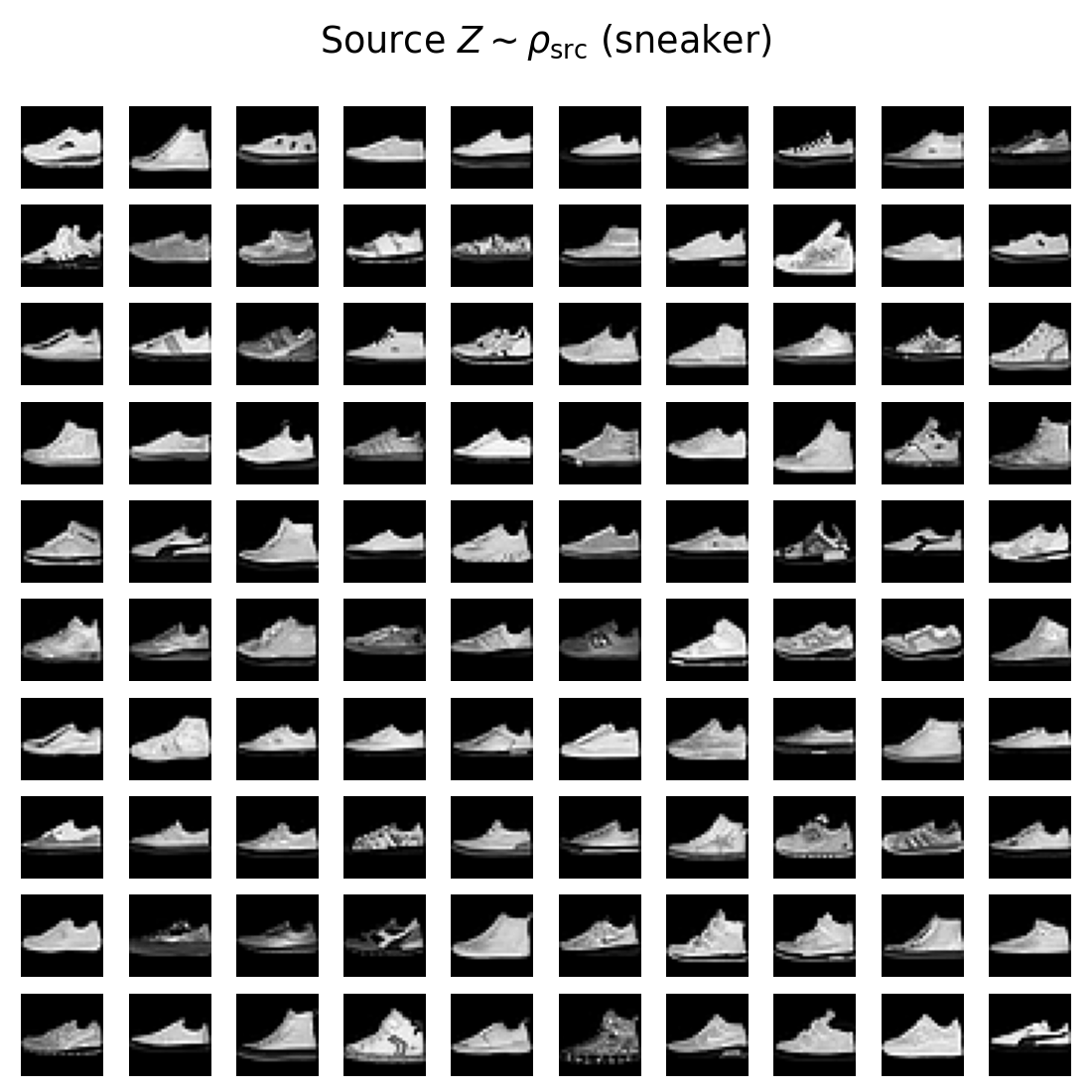}%
     }%
     \subfloat{%
         \includegraphics[width=0.495\textwidth]{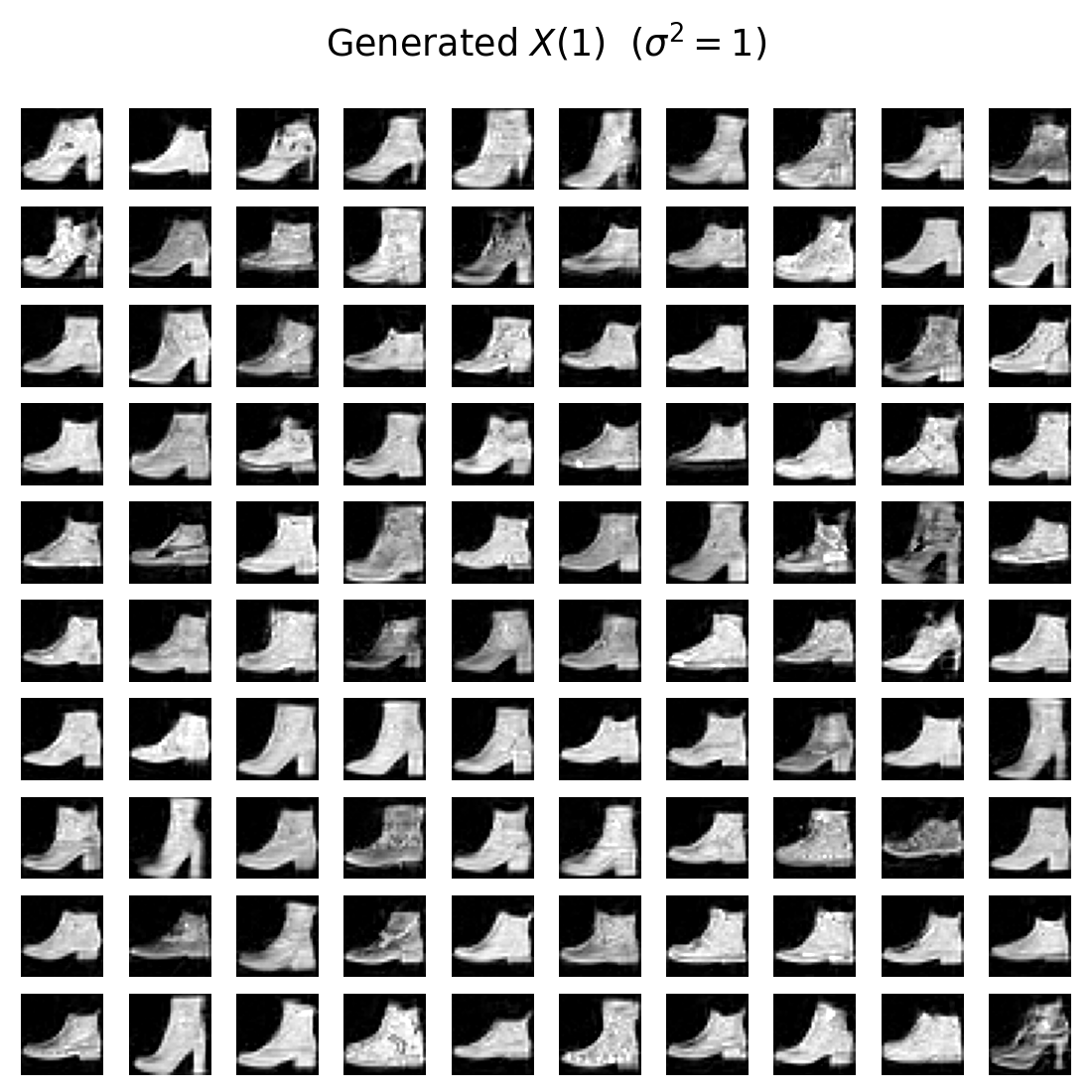}%
     }%
     \caption{
     Left panel: 100 samples from the source distribution $\rho_{\text{src}}$ (sneaker). 
        Right panel: transported samples from the left panel meant to match the target distribution $\rho_{\text{tg}}$ (ankle boot).
     }
     \label{fig:FMNIST_X1}
\end{figure}

\begin{figure}
    \centering
    \vspace{-5pt}
    \includegraphics[width=\linewidth]{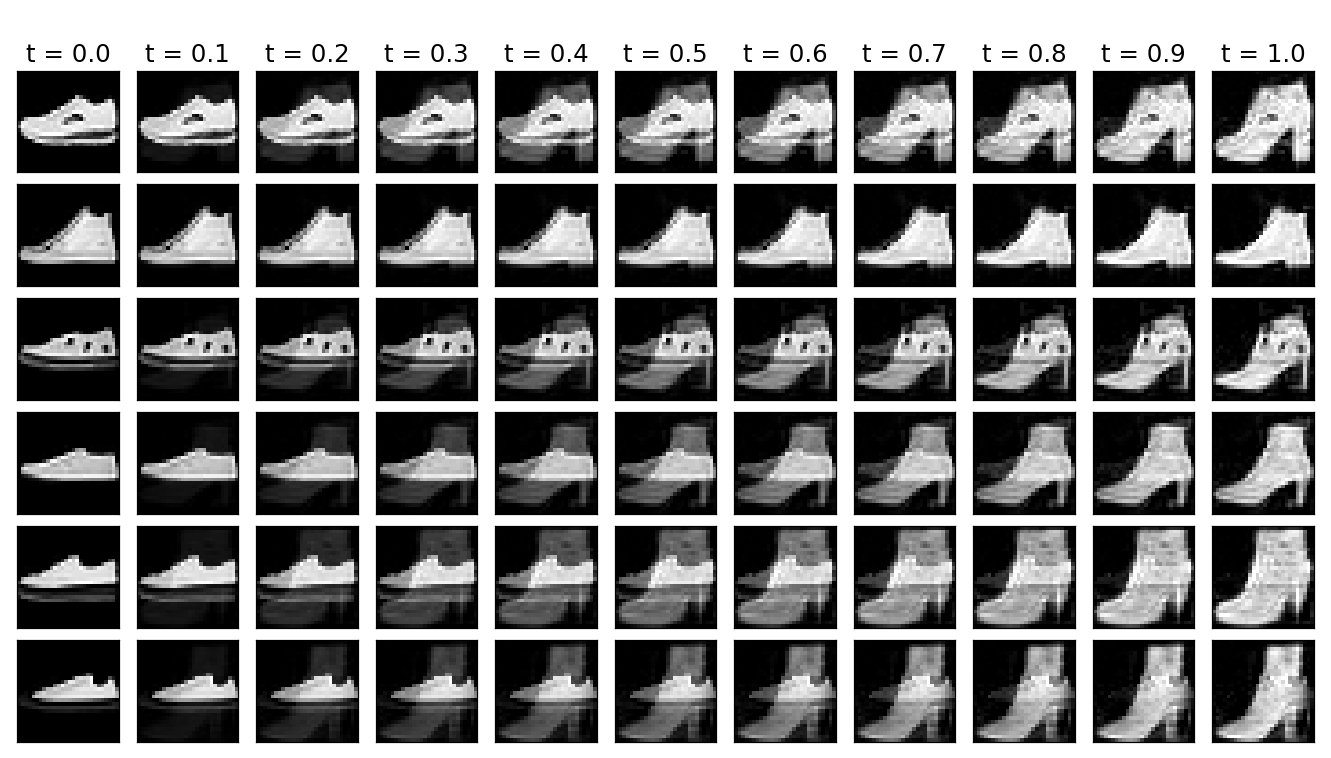}
    \captionsetup{aboveskip=0pt}
    \caption{
    Selected time slices of the transformation process $X(t)$.}
   \label{fig:FMNIST_Xt}
\end{figure}

Figure \ref{fig:FMNIST_X1} reports the generated ankle boot samples obtained from sneaker samples, while Figure \ref{fig:FMNIST_Xt} demonstrates the transformation  induced by the learned velocity field $v_{\theta^{*}}$.

\subsection{Comparison with optimal transport}\label{subsec:OTcomparison}
\begin{figure}
    \centering
    \vspace{-5pt}
    \includegraphics[width=\linewidth]{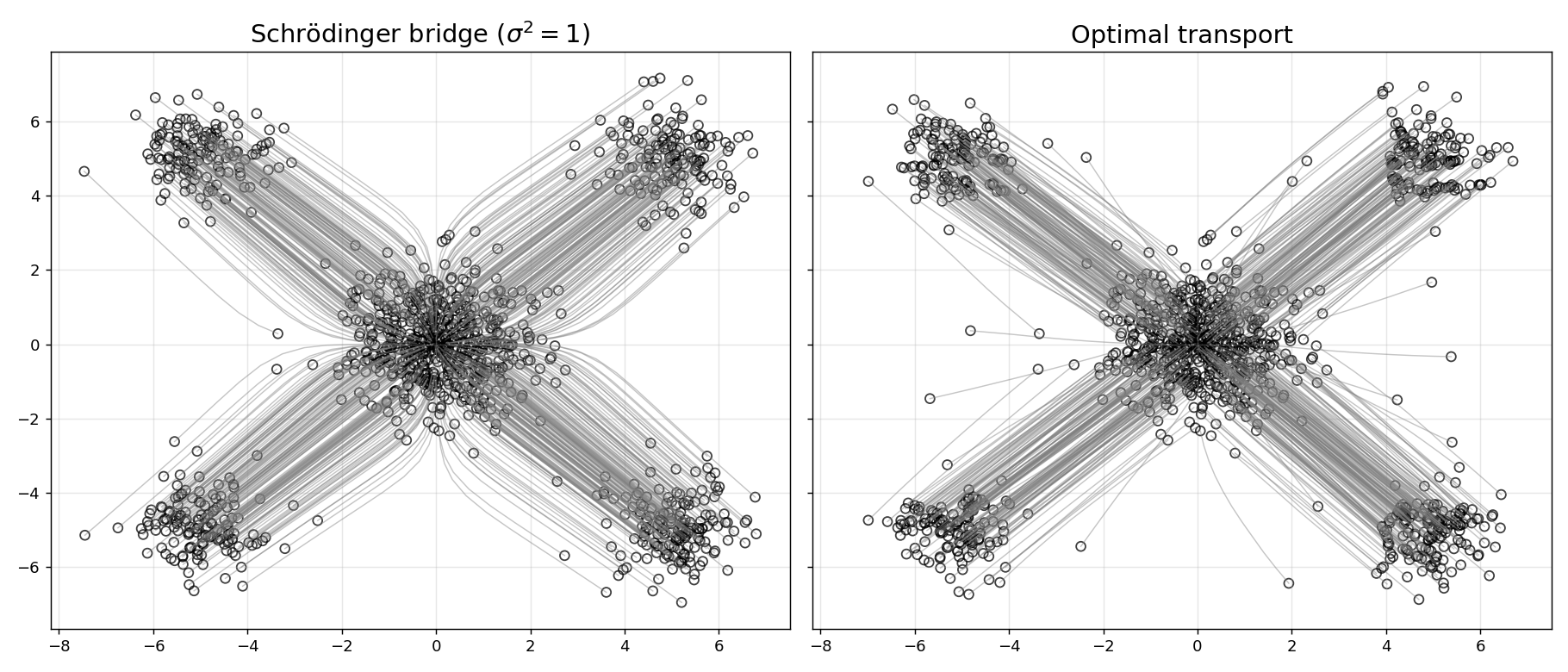}
    \captionsetup{aboveskip=0pt}
    \caption{ Left panel: learned SB trajectory copied from Figure \ref{fig:MixG_Xt}. 
        Right panel: learned OT ($\sigma^2=0$) trajectory.}
   \label{fig:SBvsOT}
\end{figure}

It is worth noting that in the case $\sigma^2=0$, the SB problem reduces to Benamou and Brenier's formulation for optimal transport (OT) \cite{benamou2000computational}:
\begin{subequations}
\begin{gather*}
    \inf_{v}\quad\E \int_0^1 \frac{1}{2} 
    \left\lVert v\bigl(X(t),t\bigl) \right\rVert^2 \dd t,\\
    \dd X(t) = v\left(X(t),t\right) \dd t,\\
         X(0)\sim \rho_{\text{src}},\, 
         X(1)\sim \rho_{\text{tg}},
\end{gather*}
\end{subequations}
which can also be solved using similar min-max neural algorithm as proposed in this work, while no score approximation is needed in the absence of the diffusion term.
However, the practical computation eventually benefits from considering the noise as the score regularization in terms of better quality of generated samples. 
To see this, in Figure \ref{fig:SBvsOT}, we compare the learned trajectory of Gaussian mixture example (Section \ref{sec:MixG2D}) in this work with that in \cite{chen2026neural}, which learns parameterized OT trajectory from source and target samples.
Starting from the identical set of source samples, we see in the right panel, the learned OT trajectories are straight lines with some points not well categorized into one of the 4 Gaussian components and fall between modes, while the learned SB trajectories in the left panel do not exhibit such pattern.

\subsection{Comments on the numerical results}\label{subsec:CommentsNumerics} 
The effectiveness of our approach is supported by the alignment of the sample cloud generated by our learned current velocity $v_\theta$ with the density contour obtained by solving the dual PDE system in \cite{caluya2021wasserstein} and \cite{11239424}, as shown in the first row of Figure \ref{fig:MixG_Xt}. 
Also, the convergence of the critic value in Figure \ref{fig:MixG_Xt_loss} agrees on the match of the terminal distribution with the target distribution. 

A notable feature of the numerical results is that the transformation is executed along a low-cost trajectory while preserving the essential structure of the source samples. 
For example, as shown in the fifth row of Figure \ref{fig:MNIST_Xt}, the ``C"-shaped stroke of the digit (9) is preserved throughout the transformation, even as the sample evolves into a digit (6). 
Such pattern is also be observed from an extra numerical experiment of human face images provided in Appendix~\ref{app:CelebA}.
Contrasting different source distributions in the two MNIST experiments Section \ref{test:MNIST} and Section \ref{test:MNIST_Gau}, the evolution starting from standard Gaussian seems to flow toward target in a path with noise (Figure \ref{fig:MNISTG_Xt}) while the other one starting without noise does not show such pattern (Figure \ref{fig:MNIST_Xt}).
We observe the quality of target sample generated under Schr\"odinger bridge model outperform that generated under OT model, where the latter is adopted by many flow matching models.

\section{Concluding Remarks}\label{sec:ConcludingRemarks}
We propose a novel neural solver for the Schr\"{o}dinger bridge problem. The solver builds on a saddle point formulation involving the Yasue action functional, and allows for learning an optimal current velocity governing the sample evolution via an ordinary -- not stochastic -- differential equation. The proposed solver relies only on the availability of the endpoint samples, and accounts for the score evolution along the trajectory. We explain the proposed formulation and detail its implementation via standard actor-critic architecture. We demonstrate the efficacy of the proposed solver using several numerical experiments involving real and synthetic datasets.

%
\section*{Conflict of interest}
The authors declare that they have no conflict of interest.

\medskip

\bibliography{references.bib}
\bibliographystyle{plain}

\appendix

\section{Proof of Proposition \ref{pr:consistent}}\label{pf:consistent}
We show that Problem \eqref{sch_fish_SDE} and Problem \eqref{pb:saddle} have the same optimal value:
\begin{equation*}
    \inf_{v}\sup_{\phi:  C_b(\Rmb^d)} f_1(v)+ f_2(v,\phi) = 
    \inf_{v:X_v(1)\sim\rho_{\text{tg}}} f_1(v).
\end{equation*}
Here, 
    \begin{align*}
        X_v(t):& = Z + \int_0^t v\left(X_v(t),t\right)\dd t, \\
    f_1(v): & =
 \int_0^1\mathbb{E}
   \left[ \frac{1}{2} \lVert v\left(X_v(t),t\right)\rVert^2 + \frac{\sigma^4}{8}\left\lVert\nabla\log\rho_v\left(X_v(t),t\right)\right\rVert^2\right]
    \dd t,\\
    f_2(v,\phi): & = \mathbb{E}\left[\phi\left(X_v(1)\right)\right]-
    \mathbb{E}\left[\phi(Y)\right],
    \end{align*}
where $\rho_v(\cdot,t)$ denotes the density of $X_v(t)$.
All infima are taken over admissible velocity fields for which the trajectories and the objective are well defined. 

If $X_v(1)\sim\rho_{\text{tg}}$, then
$f_2(v,\phi)=0$ for every $\phi\in C_b(\Rmb^d)$.
Otherwise, since bounded continuous functions distinguish
probability measures, there exists $\phi_0\in C_b(\Rmb^d)$
such that $f_2(v,\phi_0)\ne0$.
Replacing $\phi_0$ by $-\phi_0$ if necessary, we may assume
$f_2(v,\phi_0)>0$. By linearity,
\[
    f_2(v,\lambda\phi_0)
    =\lambda f_2(v,\phi_0)\longrightarrow+\infty
    \qquad\text{as }\lambda\to+\infty.
\]
Consequently,
\[
    \sup_{\phi\in C_b(\Rmb^d)} f_2(v,\phi)
    =
    \begin{cases}
        0, & X_v(1)\sim\rho_{\text{tg}},\\
        +\infty, & \text{otherwise}.
    \end{cases}
\]
Taking the infimum over $v$ proves the claimed equality.

Moreover, suppose $(v^*,\phi^*)$ is a saddle point of
$L(v,\phi):=f_1(v)+f_2(v,\phi)$, with $f_1(v^*)<\infty$.
The saddle inequalities give
\[
    L(v^*,\phi)\le L(v^*,\phi^*)\le L(v,\phi^*)
    \qquad\text{for all admissible }v
    \text{ and }\phi\in C_b(\Rmb^d).
\]
The first inequality and the characterization above imply
$X_{v^*}(1)\sim\rho_{\text{tg}}$: 
Indeed, tet $\phi$ be any function in $C_b(\Rmb^d)$. 
Then both $\phi^*+\phi$ and $\phi^*-\phi$ are in $C_b(\Rmb^d)$. 
Specifically, since $\phi^*+\phi$ is not the optimizer $\phi^*$ in general, 
\begin{equation*}
    f_1(v^*)+ f_2(v^*,\phi^*+\phi)
    \le f_1(v^*)+ f_2(v^*,\phi^*),
\end{equation*}
which further implies
\begin{equation*}
    \mathbb{E} \Bigr[  \phi\left(X_{v^*} (1)\right)\Bigr]-
    \mathbb{E}\left[\phi(Y)\right]
    \le 0,
\end{equation*}
Apply the same argument as above to $\phi^*-\phi$ to obtain 
\begin{equation*}
     \mathbb{E} \Bigr[ \phi\left(X_{v^*}(1)\right)\Bigr]-
    \mathbb{E}\left[\phi(Y)\right]
    \ge 0.
\end{equation*}
As a result, for all  $\phi\in C_b(\Rmb^d)$, we have
\begin{equation*}
    f_2(v^*,\phi) = \mathbb{E} \Bigr[  \phi\left(X_{v^*}(1)\right)\Bigr]-
    \mathbb{E}\left[\phi(Y)\right] = 0,
\end{equation*}
which further implies the feasibility of $v^*$, i.e., $X_{v^*}(1)\sim\rho_{\text{tg}}$ and 
\begin{equation}\label{ineq:ge}
    f_1(v^*)\ge \inf_{v:X_v(1)\sim\rho_{\text{tg}}}f_1(v).
\end{equation}

For any feasible $v$, both
$f_2(v^*,\phi^*)$ and $f_2(v,\phi^*)$ vanish.
The second inequality therefore yields
\[
    f_1(v^*)\le f_1(v).
\]
Thus $v^*$ solves the constrained problem, and we have 
\begin{equation}\label{ineq:le}
      f_1(v^*)\le \inf_{v:X_v(1)\sim\rho_{\text{tg}}}f_1(v).
\end{equation}
We thus conclude the proof by combining \eqref{ineq:ge} and \eqref{ineq:le}.\hfill\qedsymbol

\section{Proof of Proposition \ref{pr:score_ODE}} \label{pf:score_evo}
We derive the evolution of the score $s=\nabla \log \rho$ along trajectories associated with  the continuity equation 
\begin{equation}\label{appdx:eqn_conti}
    \partial_t \rho  + 
\nabla\cdot\left(v\rho\right)=0,
\end{equation}
whose characteristic trajectories satisfy 
\begin{equation}\label{appdx:eqn_transport}
    \dd x_t = v(x_t,t)\, \dd t.
\end{equation}
Assume that $\rho$ and $v$ are sufficiently smooth and that $\rho>0$ in the region under consideration. Throughout, we use the Jacobian convention 
$$
(\nabla v)_{ij}=\frac{\partial v_i}{\partial x_j}.
$$
Expanding the divergence in  \eqref{appdx:eqn_conti} gives  
\begin{equation*}
    \partial_t\rho(x,t) + (v\cdot\nabla\rho)(x,t)=  -(\rho\nabla\cdot v )(x,t).
\end{equation*}
Dividing by $\rho$, we obtain 
\begin{equation*}
     \partial_t \log\rho(x,t) +  (v\cdot\nabla\log\rho)(x,t) =  -(\nabla\cdot v )(x,t).
\end{equation*}
Taking the spatial gradient and setting $s=\nabla\log\rho$ yields 
\begin{equation*}
     \partial_t s(x,t) + (\nabla s(x,t))^T v(x,t)
      =  - (\nabla v (x,t))^Ts(x,t) -(\nabla\nabla\cdot v )(x,t), 
\end{equation*}
where we have used
$$
\nabla(v\cdot s)=(\nabla v)^Ts+(\nabla s)^Tv
$$
Plugging $x_t$ into the spacial argument of the score $s$, then by \eqref{appdx:eqn_transport} we finally arrive at the total derivative of $s(x_t,t)$ with respect to $t$ so that 
\begin{equation*}
     \frac{\dd s(x_t,t)}{\dd t} = -   (\nabla v (x_t,t))^Ts(x_t,t)  -(\nabla\nabla\cdot v )(x_t,t), 
\end{equation*}
where we have used the symmetry of $\nabla s=\nabla^2\log\rho$.\hfill\qedsymbol

\section{An additional numerical experiment with CelebA dataset} \label{app:CelebA}
We apply our method to the CelebA dataset~\cite{liu2015faceattributes},
which contains over $200,000$ face images of size $178 \times 218$ pixels with
three color channels, each annotated with 40 binary attributes. Using the \texttt{Smiling} attribute, we define non-smiling faces as the
source distribution $\rho_{\text{src}}$ and smiling faces as the target distribution $\rho_{\text{tg}}$.

The aligned images are center-cropped to $140\times140$ pixels and resized to
$64\times64$, yielding vectors in
$\R^{64\times64\times3} = \R^{12,288}$.
To reduce the computational cost, we construct the bridge in a latent space,
following established approaches for image data~\cite{rombach2022high,korotin2024light,vahdat2021score}. Specifically, the pretrained variational autoencoder (VAE) of Stable
Diffusion~\cite{rombach2022high} maps each image to a latent representation of dimension
$d = 4\times8\times8 = 256$, in which we train  the networks $v_\theta$ and $\phi_\omega$.
The encoder and decoder remain frozen throughout training; the decoder is used only to visualize the latent trajectories $X(t)$. As in the preceding experiments, the initial score $s\left(X(0),0\right)$ is learned beforehand by denoising score matching, since the source distribution $\rho_{\text{src}}$ is available only through samples.

\begin{figure}[htbp]
     \centering
     \subfloat{%
         \includegraphics[width=0.495\textwidth]{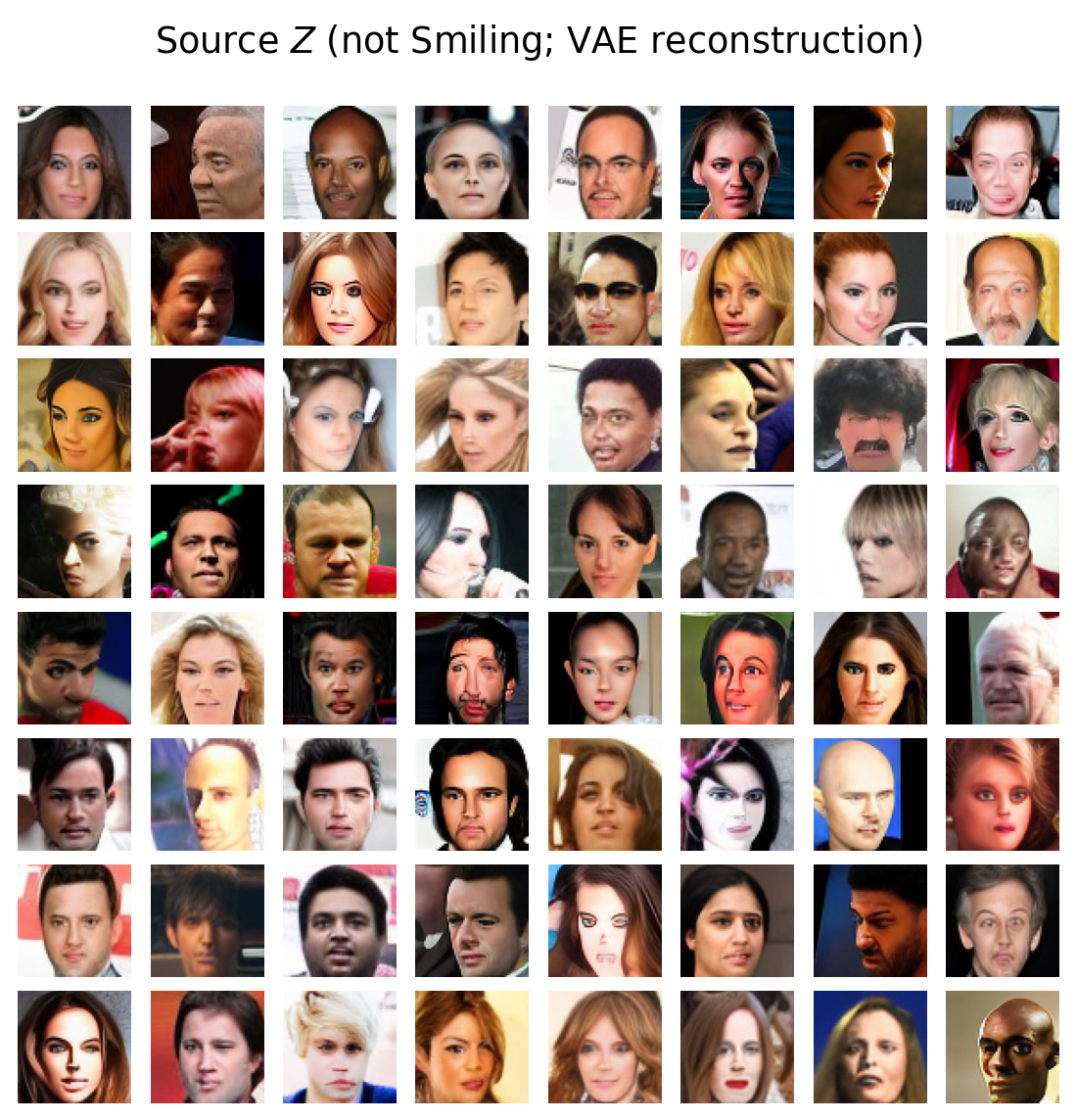}%
     }%
     \subfloat{%
         \includegraphics[width=0.495\textwidth]{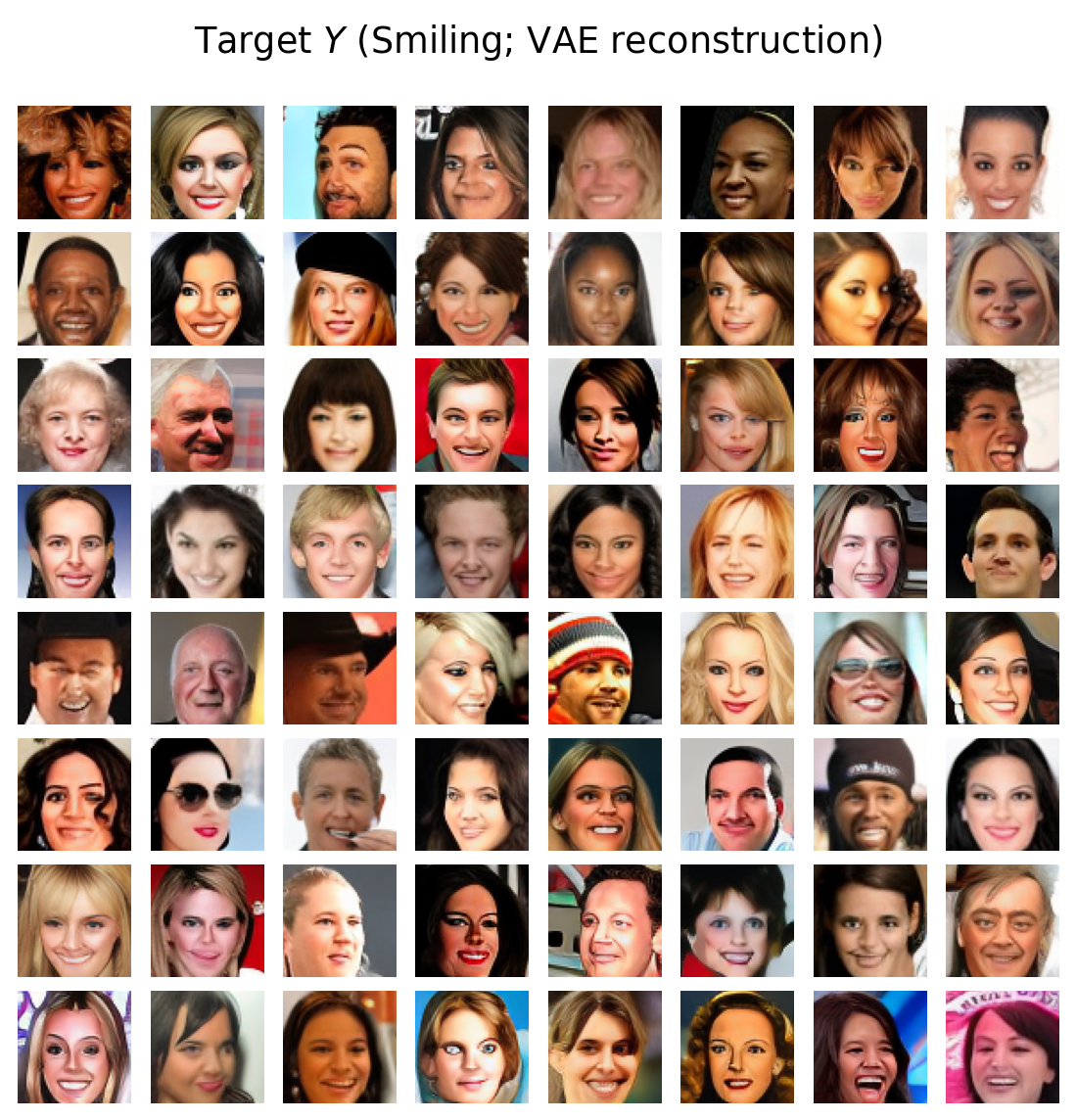}%
     }%
     \caption{
     Left panel: 64 samples from the dataset labeled \texttt{non-smiling} ($\rho_{\text{src}}$). 
        Right panel: 64 samples from the dataset labeled \texttt{smiling} ($\rho_{\text{tg}}$).
     }
     \label{fig:CelebA_Samples}
\end{figure}

\begin{figure}
    \centering
\vspace{-4pt}
\includegraphics[width=\linewidth]{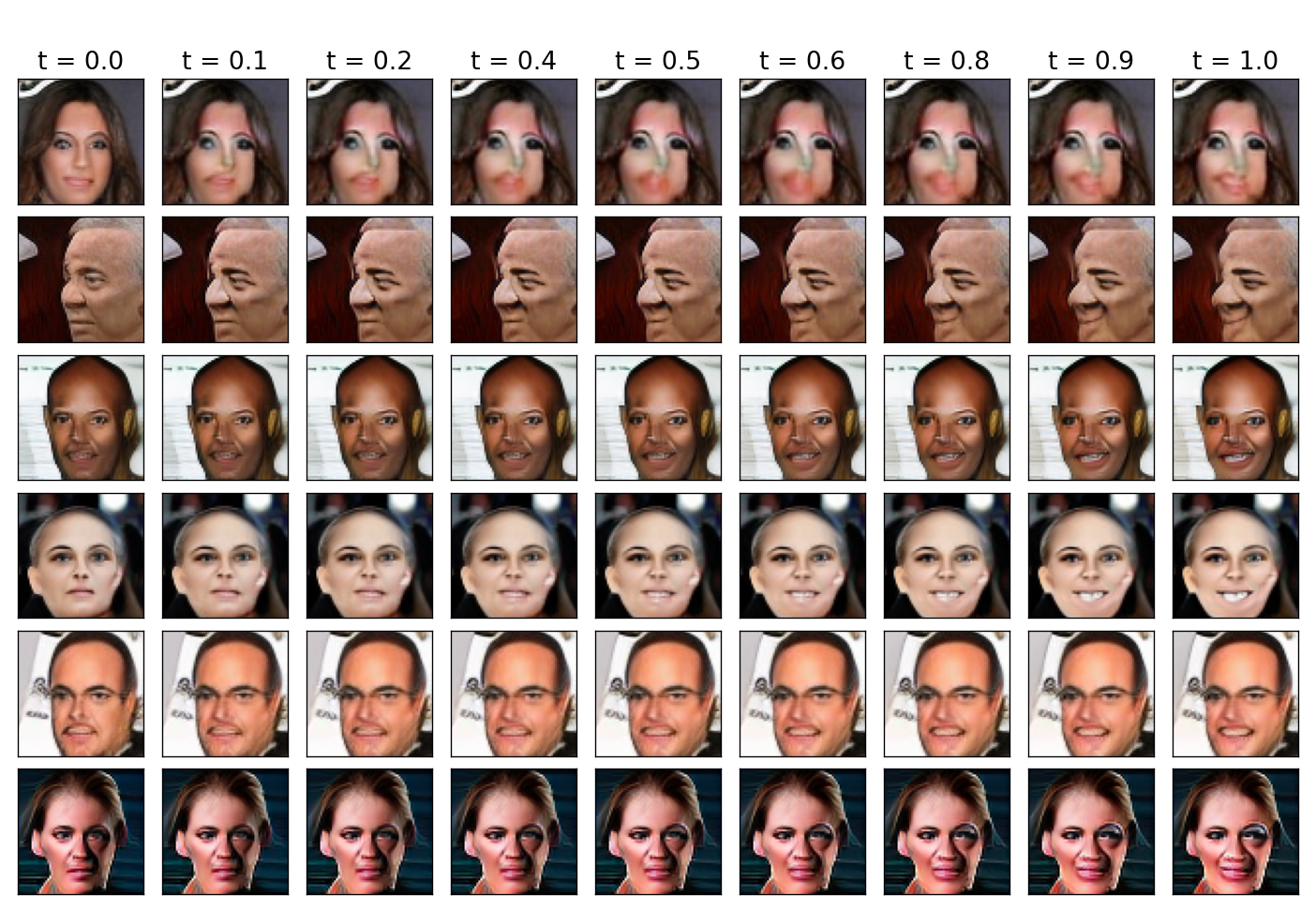}
    \captionsetup{aboveskip=-1pt} 
    \caption{Selected time slices of the transformation process $X(t)$.
       }
    \label{fig:CelebA_Xt}
\end{figure}

The source and target samples are unpaired: no correspondence between individual faces is provided, and their coupling is determined by minimizing the objective. Since smiling is
largely independent of identity, pose, and background, we expect a successful
bridge to change the expression while leaving these features unchanged. Figure~\ref{fig:CelebA_Samples} shows unpaired samples from both distributions,
displayed as VAE reconstructions (encoded and then decoded). Figure~\ref{fig:CelebA_Xt} shows the learned transformation from non-smiling
($\rho_{\text{src}}$) to smiling ($\rho_{\text{tg}}$) faces.
 
\section{Network architectures and training hyperparameters}\label{appdx:details_exp_hyper}

Since $\phi\in C_b(\mathbb{R}^d)$, we choose the critic network to be of the form
$$\phi_\omega(x)=\phi_{\max}\tanh\left(L_\phi h_\omega(x)/\phi_{\max}\right),$$ with $h_\omega$ 1-Lipschitz via spectral normalization \cite{miyato2018spectral}, so
$\|\phi_\omega\|_\infty\le\phi_{\max}$ and Lipschitz constant at most $L_\phi$.
In Table~\ref{tab:hyperparameters} and Table~\ref{tab:conv_critic}, we provide the network architectures and the associated hyperparameters used in our experiments.


\begin{table}[htbp]
\centering
\caption{Network architectures and training hyperparameters.}
\label{tab:hyperparameters}
\footnotesize
\setlength{\tabcolsep}{6pt}
\begin{tabular}{@{}lccccc@{}}
\toprule
 & Mixed & MNIST & MNIST & Fashion- & CelebA \\
 & Gaussian & $9\to6$ & $\mathcal N(0,I)\to6$ & MNIST & Smiling \\
\midrule
\multicolumn{6}{@{}l}{$v_\theta$} \\
Network type  & \multicolumn{5}{c}{multilayer perceptron} \\
Layers $\times$ width & $5\times128$ & $5\times512$ & $5\times512$ & $5\times512$ & $5\times1024$ \\
Activation         & \multicolumn{5}{c}{Sigmoid Linear Unit (SiLU)} \\
Learning rate      & $5\times10^{-5}$ & $5\times10^{-5}$ & $5\times10^{-5}$ & $10^{-4}$ & $5\times10^{-5}$ \\
\midrule
\multicolumn{6}{@{}l}{$\phi_\omega$} \\
Network type               &  \multicolumn{4}{c}{multilayer perceptron}  & convolutional$^{\dagger}$ \\
Layers $\times$ width & $3\times128$ & $3\times512$ & $3\times512$ & $3\times512$ & -- \\
Activation         & \multicolumn{5}{c}{LeakyReLU (slope $0.2$)} \\
Learning rate      & $5\times10^{-5}$ & $5\times10^{-5}$ & $8\times10^{-5}$ & $10^{-4}$ & $5\times10^{-5}$ \\
Bound $\phi_{\max}$ & 180.0 & 6797.0 & 6797.0 & 394.7 & 48.6 \\
Lipschitz $L_\phi$ & 21.0 & 129.0 & 129.0 & 31.1 & 10.9 \\
\midrule
\multicolumn{6}{@{}l}{\textit{Training} $v_\theta$ and $\phi_\omega$} \\
Batch size $n$     & 2000 & 200 & 200 & 200 & 256 \\
Time steps $N$     & 10 & 64 & 64 & 64 & 64 \\
\midrule
\multicolumn{6}{@{}l}{\textit{Pre-trained source score $s_\psi$}} \\
Layers $\times$ width & -- & $4\times512$ & -- & $4\times512$ & $4\times1024$ \\
Activation    & \multicolumn{5}{c}{LeakyReLU (slope $0.2$)} \\
Learning rate      & -- & $10^{-4}$ & -- & $10^{-4}$ & $3\times10^{-4}$ \\
\bottomrule
\end{tabular}

\smallskip
\parbox{\linewidth}{\footnotesize
$^{\dagger}$See Table~\ref{tab:conv_critic}.}
\end{table}

\begin{table}[htbp]
\centering
\caption{Convolutional critic for CelebA, acting on the latent $z\in\mathbb{R}^{4\times8\times8}$.}
\label{tab:conv_critic}
\small
\begin{tabular}{@{}lccc@{}}
\toprule
Layer & Kernel, stride, padding & Output shape & Activation \\
\midrule
Input $z$          & --               & $4\times8\times8$   & -- \\
Conv               & $3\times3$, 1, 1 & $64\times8\times8$  & LeakyReLU (0.2) \\
Conv               & $4\times4$, 2, 1 & $128\times4\times4$ & LeakyReLU (0.2) \\
Conv               & $4\times4$, 2, 1 & $256\times2\times2$ & LeakyReLU (0.2) \\
Flatten + linear   & --               & $1$ & -- \\
\bottomrule
\end{tabular}
\end{table}

\end{document}